\documentclass[12pt]{article}
\usepackage{amssymb}
\usepackage{amsmath,amsthm}
\usepackage{amsfonts}
\usepackage{graphicx}
\usepackage{graphics}

\numberwithin{equation}{section}
\numberwithin{footnote}{section}

\newtheorem{cor}{Corollary}[section]

\newcommand{\bp}{\begin{proof}}
\newcommand{\ep}{\end{proof}}
\newcommand{\be}{\begin{equation}}
\newcommand{\ee}{\end{equation}}
\newcommand{\bes}{\begin{equation*}}
\newcommand{\ees}{\end{equation*}}
\newcommand{\mbf}{\mathbf}
\newcommand{\bsy}{\boldsymbol}

\begin{document}
\date{\small\textsl{\today}}
\title{
Radial boundary elements method, \\ a new approach on using radial basis functions \\ to solve  partial differential equations, efficiently
% and its implementation for advection diffusion equation %removing singularity of integrals by an \\ optimal selection of boundary source points 
%On optimal boundary source points of boundary elements method and its implementation for \\ radial BEM
}
\author{Hossein Hosseinzadeh
\begin{footnote}{Corresponding author.\newline {\em  E-mail addresses:}
\newline h{\_}hosseinzadeh@aut.ac.ir  ,  hosseinzadeh@pgu.ac.ir (H. Hosseinzadeh). \\
z.sedaqatjoo@aut.ac.ir , zeinab.sedaghatjoo@gmail.com (Z. Sedaghatjoo). 
}$\vspace{.2cm} $
\end{footnote} 
, Zeinab Sedaghatjoo
\\
\small{\em  Department of Mathematics, Persian Gulf University, Bushehr, Iran.}
\vspace{-1mm}} \maketitle
%-------------------------------------------------------------
\vspace{.9cm}
%%%%%%%%%%%%%%%%%%%%%%%%%%%%%%%%%%%%%%%%%%%%%%%%%%%%%%%%%%%%%%%%%%%%%%%%%%%%%%%%%%%%%%%%%%%%%%%%%%%%%%%%%%%%%%%%%%%%%

\begin{abstract}
Conventionally, piecewise polynomials have been used in the boundary elements method (BEM) to approximate unknown boundary values.
Since infinitely smooth radial basis functions (RBFs) are more stable and accurate than the polynomials for high dimensional domains, the unknown values are approximated by the RBFs in this paper.
Therefore, a new formulation of BEM, called radial BEM, is obtained.
To calculate singular boundary integrals of the new method, we propose a new distribution for boundary source points that removes singularity from the integrals. Therefore, the boundary integrals are calculated precisely by the standard Gaussian quadrature rule (GQR) with $n=16$ quadrature nodes.
Several numerical examples are presented to check the efficiency of the radial BEM versus standard BEM and RBF collocation method for solving partial differential equations (PDEs). Analytical and numerical studies presented in this paper admit the radial BEM as a perfect version of BEM which will enrich the contribution of BEM and RBFs in solving PDEs, impressively. %
\vspace{.5cm}
\end{abstract}

%%%%%%%%%%%%%%%%%%%%%%%%%%%%%%%%%%%%%%%%%%%%%%%%%%%%
%%%%%%%%%%%%%%%%%%%%%%%%%%%%%%%%%%%%%%%%%%%%%%%%%%%%
%\newpage
\section{Introduction}\label{Sec1}
The Boundary Element Method (BEM) is a numerical method that has been utilized in many branches of science and industry \cite{brebbia2017birth, sedaghatjoo2018numerical, li2023isogeometric}. The main advantage of BEM is transforming a boundary value problem (BVP) into a boundary integral equation (BIE) using Green's identities. The BEM reduces computational dimension of the problem by one \cite{Katsi2002bou, hosseinzadeh2013boundary}. This reduction is done by the use of some special radial functions, called fundamental solutions. Fundamental solutions are unbounded at source points. The unboundedness yields singular boundary integrals which reduce the accuracy of BEM if they are neglected \cite{Katsi2002bou, dehghan2011improvement}. The integrals can not be evaluated by classical techniques, accurately. Then several techniques have been proposed to deal with this drawback. For instance, analytical techniques \cite{dehghan2012calculation, bin2020boundary, hou2019three}, semi-analytical techniques \cite{han2018semi, hosseinzadeh2014simple}, transformation techniques \cite{gu2016extended, tan2020bem, gu2018investigation}, element subdivision techniques \cite{gao2016element, zhang2019binary, assari2019numerical} and other innovative techniques \cite{gu2017general, gong2017isogeometric, sun2023study} have been proposed. 
Most of these techniques are restricted to straight boundary elements, Laplace equation and polynomial interpolation. Therefore, a numerical technique applicable to curved boundary elements, wide range of equations and arbitrary boundary values approximation is still a concerning subject in this area. After analyzing the boundary integrals, a new technique is proposed in this paper that enjoys these benefits. The new technique is based on finding the best location of boundary source points to handle the singularity. Therefore, the singular integrals of BEM can be calculated by the standard Gaussian quadrature rule (GQR), precisely. %This is the main novelty of the paper.

The accuracy of the BEM mainly depends on the approximation power of basis functions defined on the boundary elements. The most common basis functions are polynomials of first, second and higher degree. The BEM obtained by these polynomials are named linear, quadratic and higher order BEM, respectively \cite{Katsi2002bou, hosseinzadeh2014simple}. These standard basis functions can be replaced by radial basis functions (RBFs). The RBFs are extensively utilized in the literature to solve engineering problems \cite{wendland2004scattered, Buhman2009, vsarler2007global, tayari2023investigation, narimani2023predicting, hosseinzadeh2023new}. Since RBFs are more accurate and stable than piecewise polynomials for data approximation in two and three dimensions \cite{wendland2004scattered} they are appropriate replacements for the polynomials. In \cite{assari2017meshless} RBFs are applied to solve the boundary integral equation arising from the Laplace equation with Robin boundary conditions. And in \cite{assari2019numerical} the integral equation is solved by local RBFs. From \cite{assari2017meshless, assari2019numerical} RBFs lead to significantly accurate results when boundary integrals are calculated, precisely. A new quadrature formula based on the interval sub-division is presented in \cite{assari2017meshless2, assari2018approximate, assari2017numerical} to compute the integrals. The error analysis presented in the papers shows the efficiency of the formula. But, the formula is complicated because of splitting the integrals. %A new technique is presented in this paper to calculate the boundary integrals, easily.
In this paper, the RBFs are applied to approximate boundary values of BEM for solving partial differential equation (PDE)
\be\label{PDE}
\nabla^2 u(\mbf{x}) + \mbf{h}.\nabla u(\mbf{x}) + \lambda u(\mbf{x}) = 0 ,
\ee
where $\mbf{x}\in\Omega\subset \mathbb{R}^2$, $\mbf{h}=[h_1 , h_2]$ is a real vector and $\lambda$ is a constant number. The new method, named radial BEM, is simple in analysis and efficient in application. Analytical and numerical studies presented in this paper show the new method is significantly more accurate than the standard BEM and more stable than the RBF collocation method when the proposed location of boundary source points is utilized.%The optimal points are found here when the standard Gaussian quadrature rules (GQR) is used for numerical integration of singular integrals.
%... are proposed in the literature to overcome this problem.  The new technique is based on finding optimal boundary source points where the (singular and non-singular) boundary integrals are calculated by Gaussian quadrature rule (GQR) with high precision. This is the main novelty of the paper. 

%... paper of pooria...

%Then the radial basis functions are applied in this paper to approximate unknown values of BEM for solving two dimensional BVPs. 
The rest of the paper is organized as follows.
A brief introduction to RBFs and their application for data approximation are presented in Section \ref{Sec12}.
The radial BEM is described in Section \ref{Sec2} for solving two dimensional PDE (\ref{PDE}).
The singular boundary integrals of BEM are studied in Section \ref{Sec3} and an error analysis is presented there to find the main source of error of numerical integration. Then optimal location of boundary source points is proposed there minimizes the error. It is shown there that the Gaussian quadrature rule with $n=16$ quadrature nodes is able to calculate the singular integrals accurately when the location of boundary source points is optimized. %The optimal source points is obtained in this section for the GQR with $n=8$ and $16$ quadrature points.
Numerical experiments presented in Section \ref{Sec4} show the efficiency of the radial BEM versus the standard BEM and the RBF collocation method, clearly.
Consequently, the paper is completed by a brief conclusion presented in Section \ref{Sec5}.

%%%%%%%%%%%%%%%%%%%%%%%%%%%%%%%%%%%%%%%%%%%%%%%%%%%%
%%%%%%%%%%%%%%%%%%%%%%%%%%%%%%%%%%%%%%%%%%%%%%%%%%%%
\section{Radial basis functions}\label{Sec12}
Using RBFs in scattered data approximation was proposed by Hardy \cite{hardy1971multiquadric}, and they are extended to many fields of research after that \cite{wendland2004scattered, sun2023radial, tyagi2023radial, mirzaee2022application}. %several researches have been done for its expansion ...He used multiquadric RBF in order to interpolate scattered data.
From a point of review, RBFs can be divided into three groups, infinitely smooth, piecewise smooth and compact support ones. Some well-known RBFs are listed in Table \ref{tab1}. Three first RBFs in this table are infinitely smooth, two RBFs after those are piecewise smooth and two last ones are compactly supported. It is well-known that infinitely smooth RBFs yield very accurate results when the unknown function is sufficiently smooth. The accuracy of infinitely smooth RBFs mainly depends on the shape parameter, $\epsilon$, and the highest accuracy is often obtained at small shape parameters \cite{wendland2004scattered, glaubitz2023towards, chen2023selection, koupaei2018finding}.
\begin{table}[ht]
  \begin{center}
    \caption{Some well-known RBFs presented in the literature. }\label{tab1}
    \begin{tabular}{llllll} % <-- Changed to S here.
      \hline
      \hline
      name & function & lower bound for $\lambda_{min}$ &  \\
      \hline
      Gaussian (GA)  & $\exp(-\epsilon^2 r^2)$     & $h_{min}^{-d} \exp(-M_d^2/(\epsilon^2 h_{min}^2)) $ &   \\
      Multiquadric (MQ) & $\sqrt{1+\epsilon^2 r^2}$   & $(\epsilon^2 h_{min})^{-d/2+1} \exp(-M_d/(\epsilon^2 h_{min}^2)) $ & \\
      Inverse Multiquadric (IMQ) & $1/\sqrt{1+\epsilon^2 r^2}$   & $( \epsilon^2 h_{min})^{-d/2} \exp(-M_d/(\epsilon^2 h_{min}^2)) $ & \\
      \hline
      Thin plate spline (TPS)  &  $r^{2n} \ln(r)$   &  $h_{min}^{2n}$ &  \\
      Poly-harmonic Spline (PHS) & $r^{2n+1}$  & $h_{min}^{2n+1}$ & \\
      \hline
      Local $C^0$ & $(1-r)_+$               & $h_{min}$  & \\
      Local $C^2$ & $(1-r)_+^3 (1+3r)$   & $h_{min}^{5}$  &  \\
      \hline
      \hline
    \end{tabular}
  \end{center}
\end{table}
%%
%In this section, we will briefly introduce radial basis functions. %In this article, we want to calculate the distance of them to radial basis functions. where $\Phi[i,k]=\phi(r_i(x_k))$ for $k=1, 2, ... , N$, $\lambda[i]=\lambda_i$ and $b[i]=b_i$.

For numerical implementation of RBFs, let $(\mbf{p}_k , f(\mbf{p}_k) )$ are given data for $k = 1, 2, . . . , N$ where $\mbf{p}_k\in\Omega\subseteq\mathbb{R}^2$ and $f(\mbf{x}_k)\in\mathbb{R}$. There is an interpolation function $l$ satisfies
$%\begin{align}\label{eq22.1}
l(\mbf{p}_k)=f(\mbf{p}_k)  
$%\end{align}
  for $k = 1, 2, . . . , K $ where 
\begin{equation}\label{eq2.2}
l(\mbf{p}) = \sum_{k=1}^K  \gamma_k \phi_k(\mbf{p}) ,
\end{equation}
for $\mbf{p}\in\Omega$. In this equation $\gamma_k$ is a constant number and 
$
\phi_k(\mbf{p})=\phi( r_k(\mbf{p}) ) 
$
is a RBF where 
$
r_k(\mbf{p}) =\| \mbf{p} - \mbf{p}_k \| .
$
Equation (\ref{eq2.2}) concludes system of linear equations 
\begin{equation}\label{eq2.4}
\boldsymbol{\Phi} \boldsymbol{\gamma} = \mathbf{f} ,
\end{equation}
 where
\begin{align}
 \boldsymbol{\Phi}[k,j]=\phi_k(\mbf{p}_j) , 
 ~~~~~ \boldsymbol{\gamma}[k]=\gamma_k ,
 ~~~~~ \mathbf{f}[k]=f(\mbf{p}_k) ,
 \end{align}
when $\mbf{p}=\mbf{p}_k$ is inserted in the equation for $k,j=1, 2, ...,K$. 
Theoretically, the obtained system is non-singular for positive definite RBFs \ref{tab1} \cite{fasshauer2007meshfree, wendland2004scattered}. However, numerical stability of the system mainly depends on separation distance \cite{wendland2004scattered}
\be\label{hmin}
h_{min} = \frac{1}{2} \min_{k\neq j} \| \mbf{p}_k - \mbf{p}_j \| ,
\ee
and a lower bound for the smallest eigenvalue of matrix $ \boldsymbol{\Phi}$, defined as
\bes
\lambda_{min} = \min_{\boldsymbol{\gamma}\in\mathbb{R}^K} ( \boldsymbol{\gamma}^T \boldsymbol{\Phi} \boldsymbol{\gamma} ) / ( \boldsymbol{\gamma}^T \boldsymbol{\gamma} ) ,
\ees
is presented in Table \ref{tab1}. In this equation superscript $T$ shows the transpose operator. Smaller values of $h_{min}$ conclude less stable systems \cite{wendland2004scattered}. Therefore, larger values of $h_{min}$ are desirable in numerical implementation. Note that, interpolation function $l$ admits spectral convergence for infinitely smooth RBFs and it converges polynomially for piecewise smooth and compactly support RBFs \cite{wendland2004scattered}. Then, researchers mostly prefer to use infinitely smooth RBFs. Table \ref{tab1} confirms the stability of the system is a function of the dimension, $d$, and it is reduced meaningfully for higher valus of $d$. The dimension of problem is reduced by one in BEM and consequently the radial BEM, described in the next section, leads to more stable systems versus RBF collocation method wherein the RBFs are applied to solve the PDE, directly \cite{fasshauer2007meshfree, shi2023local, jiang2023stabilized}. Experiment results presented in Section \ref{Sec4} verify this claim, obviously.

%Since Gaussian RBFs are positive definite and infinitely smooth, one can use them to approximate unknown values, safely. Then we mostly use them in this paper. %The new formulation, named radial BEM, is described in the next section.

%%%%%%%%%%%%%%%%%%%%%%%%%%%%%%%%%%%%%%%%%%%%%%%%%%%%
%%%%%%%%%%%%%%%%%%%%%%%%%%%%%%%%%%%%%%%%%%%%%%%%%%%%
\section{Radial BEM}\label{Sec2}
\begin{figure}%[h]
\begin{center}
\includegraphics[width=8cm]{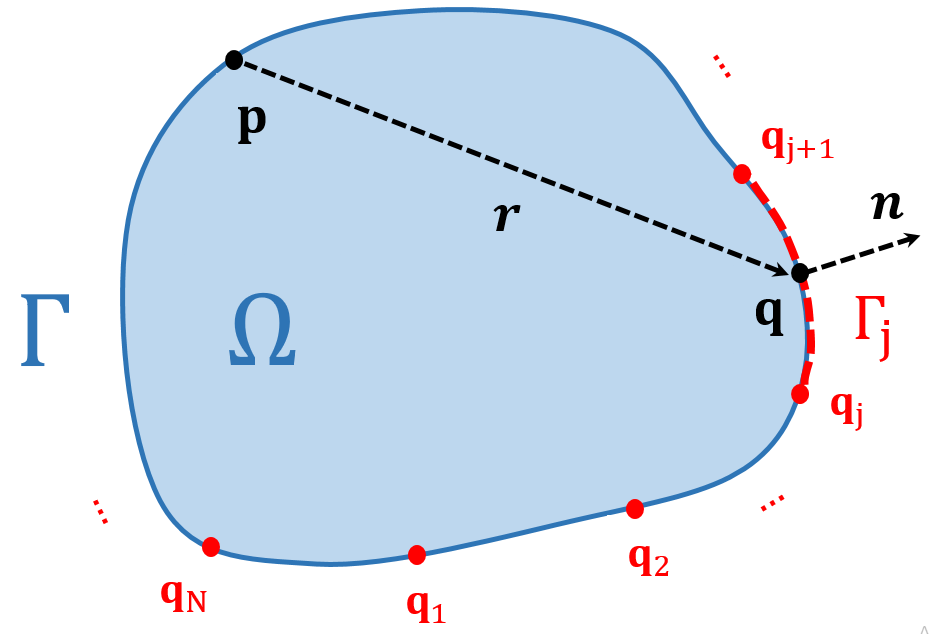}
\caption{
Computational domain $\Omega$ and its boundary, $\Gamma$. Source point, field point and distance vector are denoted by $\mbf{p}, \mbf{q}$ and $\mbf{r}$, respectively.  Boundary element  $\Gamma_j$ is obtained by connecting field points $\mbf{q}_{j+1}$ and $\mbf{q}_{j}$ to each other for $j=1, 2, ..., N$.  
}\label{fig0}
\end{center}
\end{figure}

Let $u$ be a potential function defined on the bounded computational domain, $\Omega$, with boundary $\Gamma$. The potential function may satisfy partial differential equation (\ref{PDE}) with boundary conditions
\begin{align}\label{bc}
\left\{ \begin{array}{ll}
\begin{matrix}
u(\mbf{x}) =\overline{u}(\mbf{x}) , &~~~~\hbox{for} ~~ \mbf{x}\in\Gamma_1 , \\
v (\mbf{x}) =\overline{v}(\mbf{x}) , & ~~~~ \hbox{for} ~~ \mbf{x}\in\Gamma_2 ,
\end{matrix}
\end{array} \right. 
\end{align}
where $v$ is derivative of $u$ with respect to $\mbf{n}$, i.e. $v=\partial u/ \partial \mbf{n}$, when $\mbf{n}$ is outward normal vector defined on the boundary. Also $\Gamma_1$ and $\Gamma_2$ are two parts of $\Gamma$ satisfy $\Gamma_1 \cup \Gamma_2=\Gamma$ and $\Gamma_1 \cap \Gamma_2$ is empty. PDE (\ref{PDE}) can be transformed into a boundary integral equation by the use of some special kernels.
Thanks to Divergence theorem and Green's second identity \cite{Katsi2002bou, hosseinzadeh2013boundary}, PDE (\ref{PDE}) is converted to
\be\label{eq2.1}
-c(\mbf{p}) u(\mbf{p}) = \int_{\Gamma} 
 v^*(\mbf{r}) u(\mbf{q}) d\Gamma  - \int_{\Gamma}  u^*( \mbf{r})   u(\mbf{q}) \,  \mbf{h}.\mbf{n}\, d\Gamma 
- \int_{\Gamma}  u^*( \mbf{r}) v(\mbf{q}) d\Gamma  ,
\ee
for a source point $\mbf{p}\in\mathbb{R}^2$, filed point $\mbf{q}\in\Gamma$ and $\mbf{r} =\mbf{q}-\mbf{p}$. Parameters $\mbf{p}$, $\mbf{q}$ and $\mbf{r}$ are shown in Figure \ref{fig0}, graphically.
In BIE (\ref{eq2.1}), we have $c(\mbf{p})=\alpha/2\pi$ where $\alpha$ is interior angle respect to the source point $\mbf{p}$. It is well-known that $c(\mbf{p})=1$ when $\mbf{p}$ is in $\Omega$, it equals $0.5$ when $\mbf{p}$ is located at a smooth part of the boundary and it vanishes when $\mbf{p}$ is out of the region \cite{Katsi2002bou}. 
Function $u^*$ is fundamental solution of the PDE and $v^*=\partial u^*/\partial \mbf{n}$.  The main advantage of the fundamental solution is that it satisfies
\bes
\nabla^2 u^*(\mbf{r}) - \mbf{h}.\nabla u^*(\mbf{r}) + \lambda u^*(\mbf{r}) = -\delta(\mbf{r}) ,
\ees
where $\delta$ is Dirac delta function satisfies $\delta(\mbf{r})=\infty$ for $\mbf{r}=(0,0)$ and $\delta(\mbf{r})=0$, elsewhere \cite{Katsi2002bou, ortner2015fundamental}. 
The fundamental solution of the PDE is presented in BEM literature as \cite{ortner2015fundamental, hosseinzadeh2020stability}
\be\label{fun}
u^*(\mbf{r})  = \frac{1}{2\pi} \exp( \mbf{h}.\mbf{r}/2) \, K_0(\mu \, r) ,
\ee
when $\| \mbf{h} \|^2 > 4 \lambda$ and $\mu=\sqrt{\| h \|^2 /4 - \lambda}$. Also, for Laplace equation, i.e. when $\mbf{h}=(0,0)$ and $\lambda=0$, we have
\[
u^*(\mbf{r})  = -\frac{1}{2\pi} \ln(r/\gamma) .
\] 
Note that $K_0$ is modified Bessel's function of the second kind with order zero, and $\gamma$ is a positive number bigger than or equal to the diameter of $\Omega$ guarantees stability of the BEM \cite{hosseinzadeh2020stability, sedaghatjoo2017uniqueness}. In these fundamental solutions $r$ is the absolute value of $\mbf{r}$, i.e. $r=\| \mbf{r} \|$.

%Fundamental solutions of PDEs are reported in \cite{Katsi2002bou, ortner2015fundamental, hosseinzadeh2020stability}.

To solve integral equation (\ref{eq2.1}), field points $\mbf{q}_1, \mbf{q}_2, ..., \mbf{q}_{N}, \mbf{q}_{N+1}$ are selected on the boundary of domain anticlockwise such that the last point overlaps with the first point, i.e. $\mbf{q}_{N+1}=\mbf{q}_{1}$. Then geometry of the boundary is approximated by boundary elements $\Gamma_1, \Gamma_2, ..., \Gamma_N$ obtained by connecting the field points, successively. A curved boundary element is shown in Figure \ref{fig0}. %After that, boundary values over the elements are approximated by some basis functions. Conventionally, polynomials of order less than or equal to $2$ are applied for the approximation \cite{hosseinzadeh2014simple, Katsi2002bou}. 
After the boundary discretization, BIE (\ref{eq2.1}) is reformed to
%\begin{align}\label{eq2.22}
%-c(\mbf{p}) u(\mbf{p}) = &\sum_{j=1:N}  \int_{\mbf{q}_j}^{\mbf{q}_{j+1}} 
% v^*(\mbf{r}_j(t)) u(q_j(t)) \, d\Gamma_j(t) \,  \nonumber\\
%& - \sum_{j=1:N}    \int_{\mbf{q}_j}^{\mbf{q}_{j+1}}  u^*(\mbf{r}_j(t))   u(q_j(t)) \, \mbf{h}.\mbf{n}_j(t)\,  d\Gamma_j(t) \\
%& - \sum_{j=1:N}  \int_{\mbf{q}_j}^{\mbf{q}_{j+1}}  u^*(\mbf{r}_j(t)) v(q_j(t)) \, d\Gamma_j(t)  \nonumber,
%\end{align}
\begin{align}\label{eq2.22}
-c(\mbf{p}) u(\mbf{p}) = &\sum_{j=1:N}  \int_{\mbf{q}_j}^{\mbf{q}_{j+1}} 
 v^*(\mbf{r}) u(\mbf{q}) \, d\Gamma_j \, \nonumber \\
&  - \sum_{j=1:N}    \int_{\mbf{q}_j}^{\mbf{q}_{j+1}}  u^*(\mbf{r})   u(\mbf{q}) \, \mbf{h}.\mbf{n} \,  d\Gamma_j   \\
&  - \sum_{j=1:N}  \int_{\mbf{q}_j}^{\mbf{q}_{j+1}}  u^*(\mbf{r}) v(\mbf{q}) \, d\Gamma_j  . \nonumber
\end{align}
%where  $q_j(t)\in\Gamma_j$  %is length of $j_-$th boundary element, $\Gamma_j$, i.e. $l_j=\|\mbf{q}_{j+1} -\mbf{q}_{j}\|$. And $\mbf{r}_j(t)=q_j(t)-\mbf{p}$ where $q_j$ is a linear function 
%varies from $\mbf{q}_j$ to $\mbf{q}_{j+1}$ when $t$ varies from $-1$ to $1$.  Also $\mbf{r}_j(t)=q_j(t)-\mbf{p}$,  % i.e. 
%\be\label{qj}
%q_j(t) = \frac{\mbf{q}_{j+1}+\mbf{q}_j}{2} +\frac{\mbf{q}_{j+1}-\mbf{q}_j}{2} t , \quad t\in[-1,1] .
%\ee
%and $\mbf{n}_j(t)$ is the outward normal vector perpendicular to boundary element $\Gamma_j$ at $q_j(t)$. See Figure \ref{fig0} for more understanding.
%Note that, %linear approximation (\ref{qj}) leads simplified 
%in equation (\ref{eq2.22}) boundary elements are assumed curved and $d\Gamma_j(t)$ is appeared in the formulation. However, one can omit this function by assuming straight boundary elements. However one can omit the linear approximation and assume more complicated curved boundary elements. In Figure \ref{fig_reg} the boundary elements are assumed straight and curved for square and flower-like domains, respectively.
Boundary functions $u$ and $v$ are approximated by basis functions $\phi_1, \phi_2, ..., \phi_N$ as
\be\label{uv0}
u(\mbf{p}) =  \sum_{k=1:N} \alpha_k \phi_k(\mbf{p}) ~,~~~ v(\mbf{p}) = \sum_{k=1:N} \beta_k \phi_k(\mbf{p}) ,
\ee
for $\mbf{p}\in\Gamma$. In special case, the basis functions can be radial, i.e. 
\bes
\phi_k(\mbf{p})=\phi(r_k(\mbf{p}))=\phi(\|  \mbf{p} - \mbf{p}_k \|) ~, 
\ees 
 for $k=1, 2, ..., N$. To enhance the stability of the RBFs, center points $\mbf{p}_k$ can be selected at $\Gamma_k$ such that separation distance $h_{min}$ is maximized \cite{wendland2004scattered}.
The best location of the center points is presented in Section \ref{Sec3}.
Note that equations (\ref{uv0}) can be write in vector form
\be\label{uv2}
u(\mbf{p}) = \bsy{\phi}(\mbf{p}) \, \bsy{\alpha}^T ~,~~~ v(\mbf{p}) = \bsy{\phi}(\mbf{p}) \, \bsy{\beta}^T ,
\ee
where $\bsy{\alpha}=[\alpha_1, \alpha_2, ... , \alpha_N]$, $\bsy{\beta}=[\beta_1, \beta_2, ... , \beta_N]$ and $\bsy{\phi}(\mbf{p})=[\phi_1 (\mbf{p}), \phi_2(\mbf{p}), ... , \phi_N(\mbf{p})]$. %This vector form helps us to express coefficient vectors $\bsy{\alpha}$ and $\bsy{\beta}$ as functions of $u$ and $v$, respectively when it is necessary.
%The exact position of the collocation point $\mbf{p}_j$ on boundary element $\Gamma_j$ is presented in Section \ref{Sec3}. 
%Let $\mbf{p}$ is a boundary collocation point and $u$, $v$ are approximated over the boundary by basis functions $\phi_1, \phi_2, ... , \phi_K$. 
By inserting equations (\ref{uv0}) in boundary integral equation (\ref{eq2.22}) we have
\begin{align}\label{eq2.3}
-c(\mbf{p}) \sum_{k=1:N} \alpha_k \phi_k(\mbf{p}) = &\sum_{k=1:N}\alpha_k \sum_{j=1:N} \int_{\mbf{q}_j}^{\mbf{q}_{j+1}} 
 v^*(\mbf{r}) \phi_k(\mbf{q}) \, d\Gamma_j \nonumber\\
& - \sum_{k=1:N} \alpha_k \sum_{j=1:N}   \int_{\mbf{q}_j}^{\mbf{q}_{j+1}}   u^*(\mbf{r})   \phi_k(\mbf{q}) \,  \mbf{h}.\mbf{n} \, d\Gamma_j ~ \\
& - \sum_{k=1:N} \beta_k \sum_{j=1:N}   \int_{\mbf{q}_j}^{\mbf{q}_{j+1}}  u^*(\mbf{r}) \phi_k (\mbf{q}) \, d\Gamma_j \nonumber ,
\end{align}
which it can be simplified to vector form
\begin{align}\label{eq2.4}
-c(\mbf{p}) \bsy{\phi}(\mbf{p}) \bsy{\alpha}^T=  \mbf{a}^1_{\mbf{p}} \, \bsy{\alpha}^T 
 -  \mbf{a}^2_{\mbf{p}} \, \bsy{\alpha}^T   
 -  \mbf{b}_{\mbf{p}} \, \bsy{\beta}^T  ,
\end{align}
where $k_-$th element of vectors $\mbf{a}^1_{\mbf{p}}, \mbf{a}^2_{\mbf{p}}$ and $\mbf{b}_{\mbf{p}}$ is evaluated as
\begin{align}\label{eq2.55}
\mbf{a}^1_{\mbf{p}} [k] & = \sum_{j=1:N} \int_{\mbf{q}_j}^{\mbf{q}_{j+1}} 
 v^*(\mbf{r}) \phi_k(\mbf{q}) \, d\Gamma_j  ,~~~~~~~~~ \nonumber  \\
\mbf{a}^2_{\mbf{p}}[k]  & = \sum_{j=1:N}   \int_{\mbf{q}_j}^{\mbf{q}_{j+1}}   u^*(\mbf{r})   \phi_k(\mbf{q}) \,  \mbf{h}.\mbf{n} \, d\Gamma_j   , \\
\mbf{b}_{\mbf{p}}[k]  &= \sum_{j=1:N}   \int_{\mbf{q}_j}^{\mbf{q}_{j+1}}  u^*(\mbf{r}) \phi_k (\mbf{q}) \, d\Gamma_j , ~~~~~~~~~ \nonumber 
\end{align}
for $k=1, 2, ... , N$. 
Respect to Equation (\ref{eq2.55}), generally, there are two kinds of boundary integrals 
\be\label{I12}
I_1= \int_{-1}^{1}  v^*(\mbf{r}(t)) \, g(t) \, dt ~,~~~ 
I_2= \int_{-1}^1  u^*(\mbf{r}(t))  \, g(t) \, dt~,
\ee
which have to be evaluated for obtaining the vectors. In these integrals, $g(t)$ is a smooth function defined on $[-1,1]$.
A new approach is proposed in Section \ref{Sec3} calculates the boundary integrals, easily. In that approach the boundary integrals are approximated by two weighted summations obtained from Gaussian quadrature rule. In fact %The optimal boundary source points are obtained in the new approach, and boundary integrals (\ref{I12}) are found, accurately. If the union of the quadrature nodes and weights respect to the new approach are inserted in 
two vectors $\mbf{q}_G$ and $\mbf{w}_G$ are introduced there satisfy%, respectively, we have  
%, named $\mbf{q}_G$. The Gaussian weights correspond to the Gaussian nodes are also inserted in a vector, say $\mbf{w}_G$. Then they yield
\be\label{intg}
\int_{\Gamma} f(\mbf{q}) ~d\Gamma =  
\sum_{j=1:N} \int_{\mathbf{q}_j}^{\mathbf{q}_{j+1}} f(\mbf{q}) ~d\Gamma_j
\simeq \sum_{m=1:n N} f(\mbf{q}_G[m]) \mbf{w}_G[m]= f(\mbf{q}_G) \mbf{w}_G^T ,
\ee
for a real function $f$. In this approximation, $f(\mbf{q}_G)$ is the value of $f$ at quadrature points $\mbf{q}_G$.  %Since vectors $\mbf{q}_G$ and $\mbf{w}_G$ are utilized frequently in this paper, we call them quadrature nodes and weights, respectively. 
Note that new vectors $\mbf{q}_G$ and $\mbf{w}_G$, named quadrature points and quadrature weights respectively, only depend on the boundary of the domain and they can be computed quickly when the boundary is Lipschitz. The quadrature nodes are shown in Figure \ref{fig_squ} for square and flower-like domains when $N=8$ and $n=16$. 
%where one dimensional boundary integrals can be calculated numerically by the Gaussian quadrature rule as Equation (\ref{intg}) when boundary source points are optimal. Finding the optimal center points is the main subject of Section  \ref{sec.op} and we only focus on obtaining the final system of linear equations of the new formulation in this section. 
%%%%%%%%%%%%%%%%%%%%%%%%%%%%%%%
%\subsection{The final system of linear equations}
If boundary integrals presented in Equation (\ref{eq2.55}) are approximated by quadrature rule (\ref{intg}) then Equation (\ref{eq2.55}) is converted to
\begin{align}\label{eq2.5}
\mbf{a}^1_{\mbf{p}} [k] & = (~ v^*(\mbf{q}_G-\mbf{p}) .* \phi_k(\mbf{q}_G) ~) \mbf{w}_G^T  ,~~~~~~~~~~~~~~~~~\nonumber  \\
\mbf{a}^2_{\mbf{p}} [k] & = ( ~u^*(\mbf{q}_G-\mbf{p}) .* \phi_k( \mbf{q}_G) .*  (\mathbf{h}.\mbf{n}(\mbf{q}_G)) ~) \mbf{w}_G^T  , \\
\mbf{b}_{\mbf{p}} [k] &= ( ~u^*(\mbf{q}_G-\mbf{p}) .* \phi_k( \mbf{q}_G) ~) \mbf{w}_G^T , ~~~~~~~~~~~~~~~~\nonumber 
\end{align}
where notation $.*$ shows Hadamard (or dot) product. For two matrices $H$ and $G$ of the same size, Hadamard product $H .* G$ is a matrix which is evaluated as 
\[
(H .* G)[i,j]=H[i,j] G[i,j].
\]
It should be mentioned that $\mbf{n}(\mbf{q}_G)$, presented in Equation (\ref{eq2.5}), is a matrix with two rows. The first and second rows of the matrix are devoted to the first and second components of the normal vector, $\mbf{n}$, at the quadrature points, respectively. %And $\bsy{\Phi}[:,k]$ denotes $k_-$th column of matrix $\bsy{\Phi}$. %Also $(\phi_k)_\mbf{n}(\mbf{q}_G)$ is derivative of function $\phi_k$ with respect to $\mbf{n}$ at $\mbf{q}_G$.

%\begin{align}\label{eq2.4}
%-c(\mbf{p}) \bsy{\phi}(\mbf{p}) \bsy{\alpha}^T 
%=   \mbf{H}_{\mbf{p}}   \bsy{\alpha}^T 
% -  \mbf{G}^1_{\mbf{p}}   \bsy{\alpha}^T
% - \mbf{G}^2_{\mbf{p}} \bsy{\beta}^T ,
%\end{align}
%where vectors $\mbf{H}_{\mbf{p}} , \mbf{G}^1_{\mbf{p}} , \mbf{G}^2_{\mbf{p}}$ and matrix $\bsy{\Phi}$ are evaluated as

%for $k=1, 2, ... , K$. 

%As is mentioned in section \ref{Sec3}, 
To find unknown vectors $\bsy{\alpha}$ and $\bsy{\beta}$ we need $2N$ equations. The first $N$ equations are obtained by inserting boundary source points $\mbf{p}_1, \mbf{p}_2, ... , \mbf{p}_N$ in Equation (\ref{eq2.4}) instead of $\mbf{p}$, and the last $N$ equations are obtained via boundary conditions (\ref{bc}). 
By inserting the boundary source points in Equation (\ref{eq2.4}), system of linear equations 
\begin{align}\label{eq2.6}
-\frac{1}{2} \bsy{\Phi} \bsy{\alpha}^T 
=   \mbf{A}^1  \bsy{\alpha}^T 
 -  \mbf{A}^2   \bsy{\alpha}^T
 - \mbf{B} \bsy{\beta}^T ,
\end{align}
is obtained where $k_-$th row of matrices $\bsy{\Phi}, \mbf{A}^1, \mbf{A}^2$ and $\mbf{B}$ are evaluated as
\begin{align*}
\bsy{\Phi}[k,:] =\bsy{\phi}(\mbf{p}_k) ~,~~~ \mbf{A}^1[k,:]=\mbf{a}^1_{\mbf{p}_k}~, ~~~
\mbf{A}^2[k,:]=\mbf{a}^2_{\mbf{p}_k} ~,~~~ \mbf{B}[k,:]=\mbf{b}_{\mbf{p}_k}.
\end{align*}
%when $k$ varies from $1$ to $K$. 
If we define vectors $\mbf{u}$ and $\mbf{v}$ as 
\begin{align}\label{uv}
\begin{matrix}
\mbf{u} =[u(\mbf{p}_1), u(\mbf{p}_2), & ... , u(\mbf{p}_N)] ,\\
\mbf{v} =[v(\mbf{p}_1), v(\mbf{p}_2), & ... , v(\mbf{p}_N)] ,
\end{matrix}
\end{align} 
then boundary conditions (\ref{bc}) yield
\begin{align}\label{bc2}
\left\{ \begin{array}{ll}
\begin{matrix}
\mbf{u}[k] =\bar{u}(\mbf{p}_k) ,~ & ~\hbox{for} ~ k=1, 2, ..., N_1 ,~~~~~~~~~~~~~ \\
\mbf{v}[k] =\bar{v}(\mbf{p}_k) ,~ & ~\hbox{for} ~ k=N_1+1, N_1+2, ..., N ,
\end{matrix}
\end{array} \right. 
\end{align}
when boundary points $\{ \mbf{p}_1, \mbf{p}_2, ... , \mbf{p}_{N_1} \}$ and $\{ \mbf{p}_{N_1+1}, \mbf{p}_{N_1+2}, ... , \mbf{p}_{N} \}$ are located on $\Gamma_1$ and $\Gamma_2$, respectively. 
Equation (\ref{bc2}) can be presented in matrix form as
\be\label{bc3}
\mbf{C}_1 \mbf{u} + \mbf{C}_2 \mbf{v} = \bar{\mbf{w}}^T ,
\ee
where diagonal matrices $\mbf{C}_1$ and $\mbf{C}_2$ and vector $\bar{\mbf{w}}$ are evaluated as $\mbf{C}_1[k,k]=1$ for $k=1, 2, ..., N_1$ and $0$ otherwise, $\mbf{C}_2[k,k]=1$ for $k=N_1+1, N_1+2, ..., N$ and $0$ otherwise, $\bar{\mbf{w}}[k]=\bar{u}(\mbf{p}_k)$ for $k=1, 2, ..., N_1$ and $\bar{\mbf{w}}[k]=\bar{v}(\mbf{p}_k)$ for $k=N_1+1, N_1+2, ..., N$.
Besides, Equation (\ref{uv2}) leads matrix form 
\be\label{alp}
 \mbf{u}^T = \bsy{\Phi} \bsy{\alpha}^T ~,~~~ \mbf{v}^T  = \bsy{\Phi} \bsy{\beta}^T ,
\ee
when $\mbf{p}_1, \mbf{p}_2, ... , \mbf{p}_{N}$ are inserted instead of $\mbf{p}$ in the equation.
Now, if Equation (\ref{eq2.6}) is presented as
\be\label{AB}
\mbf{A} \bsy{\alpha}^T =  \mbf{B} \bsy{\beta}^T ,
\ee
for
$
\mbf{A}=\frac{1}{2} \bsy{\Phi}  +   \mbf{A}^1   -  \mbf{A}^2 ,
$
then equations (\ref{bc3}), (\ref{alp}) and  (\ref{AB}) lead to the final system of linear equations
\be\label{final}
\begin{bmatrix} 
\mbf{A}   & -\mbf{B}     & ~\mbf{0}       & ~\mbf{0}  \\ 
\bsy{\Phi} & ~\mbf{0}    & -\mbf{I}        & ~\mbf{0}  \\ 
\mbf{0}   & ~\bsy{\Phi}  & ~\mbf{0}       & -\mbf{I}  \\ 
\mbf{0}   & ~\mbf{0}     & ~\mbf{C}_1   & ~\mbf{C}_2  \\ 
\end{bmatrix}
\begin{bmatrix} 
\bsy{\alpha}^T \\
\bsy{\beta}^T \\
\mbf{u}^T \\
\mbf{v}^T \\
\end{bmatrix}
=
\begin{bmatrix} 
\mbf{0}^T \\
\mbf{0}^T \\
\mbf{0}^T \\
\bar{\mbf{w}}^T \\
\end{bmatrix} ,
\ee
where $\mbf{0}$ and $\mbf{0}^T$ are zero matrix and vector, respectively, and $\mbf{I}$ is the identity matrix. 
After solving the final system and obtaining vectors $\bsy{\alpha}$ and $\bsy{\beta}$, from Equation (\ref{eq2.4}), the potential function $u$ can be estimated at an internal point $\mbf{p}\in\Omega$ as
\begin{align}\label{eq2.7}
u(\mbf{p})
= -  \mbf{a}^1_{\mbf{p}}  \bsy{\alpha}^T 
 +  \mbf{a}^2_{\mbf{p}} \bsy{\alpha}^T
 + \mbf{b}_{\mbf{p}} \bsy{\beta}^T ,
\end{align}
where vectors $\mbf{a}^1_{\mbf{p}}$, $\mbf{a}^2_{\mbf{p}}$ and $\mbf{b}_{\mbf{p}}$ already are defined in Equation (\ref{eq2.5}).

%%%%%%%%%%%%%%%%%%%%%%%%%%%%%%%%%%%%%%%%%%%%%%%%%%%%
%%%%%%%%%%%%%%%%%%%%%%%%%%%%%%%%%%%%%%%%%%%%%%%%%%%%
\section{Treating singularity of boundary integrals}\label{Sec3} % Gaussian quadrature rule for BEM
%Gaussian quadrature rule is accurate only if integrand $f$ in Equation (\ref{intg}) is sufficiently smooth \cite{Katsi2002bou}. Then it loses its efficiency for singular integrals (\ref{I12}). In this section we show these integrals also can be estimated by the rule, accurately if the boundary source points are selected, appropriately. 
Fundamental solution $u^*$ and its derivative respect to $\mbf{n}$, $v^*$, may be unbounded at $\mbf{r}_j=(0,0)$. So conventional numerical schemes, such as the standard Gaussian quadrature rule \cite{Katsi2002bou, arnold1991concise}, are not able to calculate the boundary integrals. Therefore, a new arrangement is suggested for boundary source points in this paper to overcome this problem. 

Let the boundary integrals be calculated by the GQR with $n$ quadrature nodes on each boundary element $\Gamma_j$ for $j=1, 2, ..., N$. Then, totally, $n \times N$ boundary points are selected to use the rule. These points are shown in Figure \ref{fig_squ} for $N=8$ and $n=16$.
The points are inserted in a vector named quadrature points and dented by $\mbf{q}_G$. The weights with respect to these points are also inserted in a vector named quadrature weights and denoted by $\mbf{w}_G$. The introduced vectors can be applied for approximation as is stated in Equation (\ref{intg}).
Note that the approximation in Equation (\ref{intg}) is only valid for smooth functions, and it will fail for unbounded ones. As was said, boundary integrals of BEM are singular at source points and consequently, the GQR loses its efficiency for them. In this section, we show the integrals can be estimated by Equation (\ref{intg}), accurately if the boundary source points are selected, appropriately.
Therefore, a set of boundary source points, $\{ \mbf{p}_1, \mbf{p}_2, ..., \mbf{p}_N \}$, is proposed here to guarantee the accuracy of Approximation (\ref{intg}) when $\mbf{p}_j\in\Gamma_j$ for $j=1, 2, ... , N$. 
Mathematically, since the GQR fail for some $s\in[-1,1]$ when $\mbf{p}_j=q_j(s)$, we are looking for the best value of $s$ satisfies
\begin{align} 
I_1 &= \int_{-1}^1     v^*(q_j(t)-q_j(s)) \, g(t) \, dt \simeq \sum_{i=1:n} v^*(q_j(t_i)-q_j(s)) \, g(t_i) w_i ~, \label{I1} \\
I_2 &= \int_{-1}^1     u^*(q_j(t)-q_j(s)) \, g(t) \, dt \simeq \sum_{i=1:n} u^*(q_j(t_i)-q_j(s)) \, g(t_i) w_i  ~,\label{I2}
\end{align}
where $t_1, t_2, ... , t_n$ and $w_1, w_2, ... , w_n$ are the standard Gaussian quadrature nodes and weights, respectively \cite{arnold1991concise}. We show in forthcoming subsections that the seeking value does not depend on functions $v^*, u^*$ and $g$, essentially, and it only depends on the number of the quadrature nodes for $n\geq 16$.
Then, the main source of error of approximations (\ref{I1}) and (\ref{I2}) is recognized in Subsection \ref{sec.error}, and after that, the optimal value of $s$ minimizes the error is obtained in Subsection \ref{sec.op}. Since the optimal value of $s$ depends on $n$, an optimal value of $n$ is suggested in Subsection \ref{sec.op2}.   

%%%%%%%%%%%%%%%%%%%%%%%%%%%%%%%%%%%%%%%%%%%%%%%%%%%%
\subsection{Error analysis} \label{sec.error}
\begin{figure}%[h]
\begin{center}
\includegraphics[width=13cm]{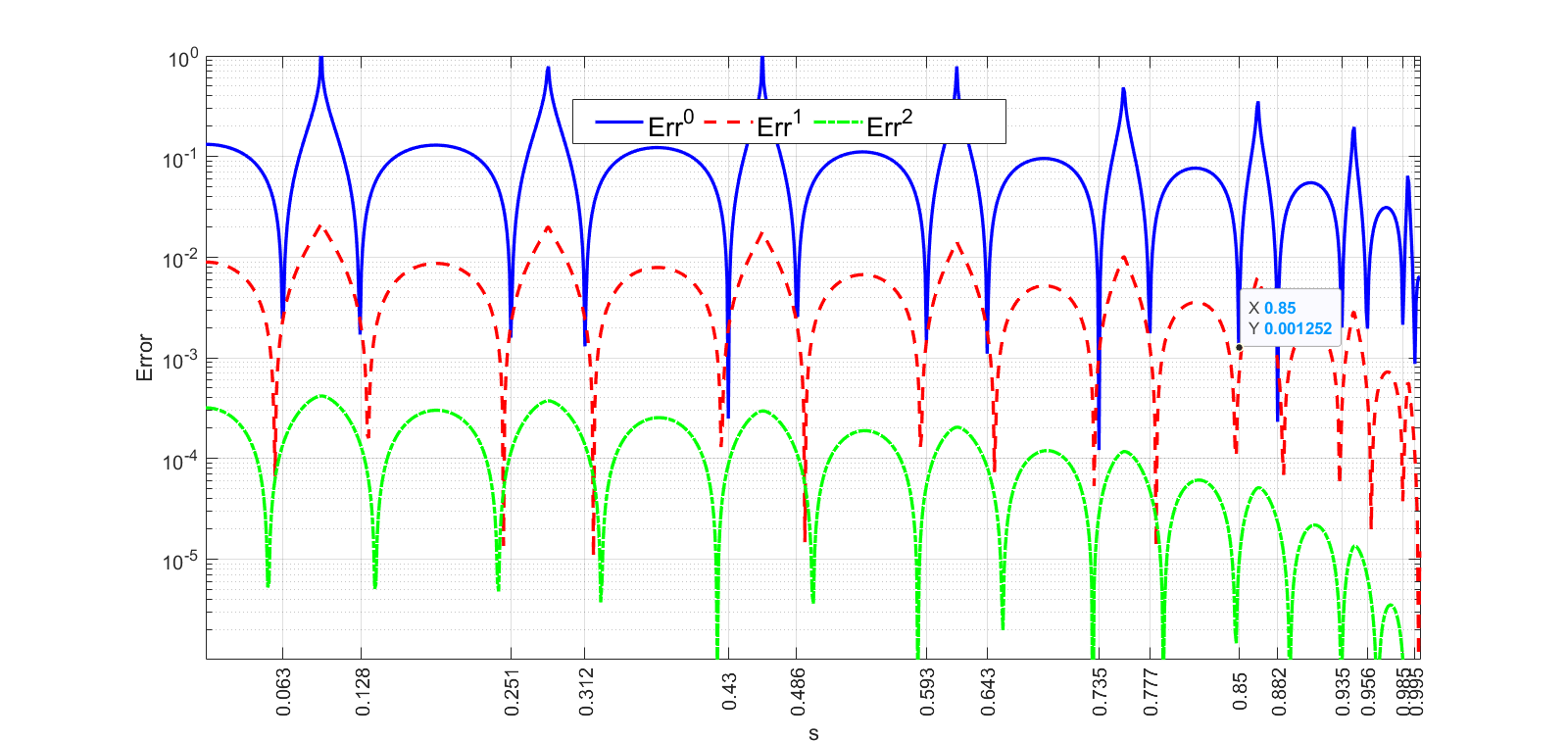}
\caption{
Graph of $Err^0$, $Err^1$ and $Err^2$ in field of $s\in[0,1]$ for $n=16$. One can see $Err^2 \leq Err^1 \leq Err^0(s)$ almost everywhere, and $Err^0(s)$ is vanished at $s= 0.063, 0.128, 0.251,$ $ 0.312, 0.430, 0.486, 0.593,$ $ 0.643, 0.735, 0.777, 0.850, 0.882, 0.935,$ $ 0.956, 0.985, 0.995$.
}\label{fig2}
\end{center}
\end{figure}

From Equation (\ref{fun}), the fundamental solution of the PDE has modified Bessel function $K_0( \mu \, r )$ in its formula. Since 
\be\label{k0}
 K_0(\mu \, r) \simeq  - \ln( r ) - \ln(\mu/2) - 0.5772 ,
\ee
when $r$ is small, it can be found the fundamental solution is not bounded at $ r=0$ \cite{hosseinzadeh2013boundary}.  In fact
\begin{align*}
u^*(\mbf{r}) & \simeq - \exp( \mbf{h}.\mbf{r}/2) \ln( r )  \simeq - \ln( r ) , \\
v^*(\mbf{r}) & \simeq - \frac{1}{2}\mbf{h}.\mbf{n} \exp( \mbf{h}.\mbf{r}/2)  \ln( r ) - \exp( \mbf{h}.\mbf{r}/2)  \frac{1}{r} \mbf{r}.\mbf{n}   \simeq -  \frac{1}{2}\mbf{h}.\mbf{n}  \ln( r ) , 
\end{align*}
for small values of $r$. 
These semi-equalities highlight the magnitude of the fundamental solution and its derivative with respect to $\mbf{n}$ increase as rate as $-\ln(r)$ when $r$ tends to zero. Then semi-equalities (\ref{I1}) and (\ref{I2}) will be accurate if and only if semi-equality
\begin{align} \label{I22}
\int_{-1}^1     \ln(\| q_j(t)-q_j(s) \|) \, g(t) \, dt \simeq \sum_{k=1:n} \ln(\| q_j(t_k)-q_j(s) \|) \, g(t_k) w_k  ~,
\end{align}
is accurate.
We know $q_j(t)$ tends to $q_j(s)$ when $t$ tends to $s$, then $q_j(t)-q_j(s)= h(t) (t-s)$ for a continuous function $h$, and approximation (\ref{I22}) will be accurate if semi-equality
\begin{align} \label{I222}
\int_{-1}^1     \ln(\| t -s \|) \, g(t) \, dt \simeq \sum_{k=1:n} \ln(\| t_k- s \|) \, g(t_k) w_k  ~,
\end{align}
is accurate.
% if collocation point $\mbf{p}$ is selected sufficiently far from field points $\mbf{q}_j$ for $j=1, 2, ..., N$. 
Thus, we are looking for that $s$ minimizes the error function
\begin{align} \label{Er}
Err(s)=\left| \int_{-1}^1     \ln(\| t-s \|) \, g(t) \, dt - \sum_{k=1:n} \ln(\| t_k-s  \|) \, g(t_k) w_k \right| ~.
\end{align}
Note that the error function mainly depends on the number of quadrature nodes, $n$.
It also  slightly depends on smooth function $g$. To omit the role of $g$ in Equation (\ref{Er}), the function can be expanded around $s$ by Taylor's expansion formula as
\bes
g(t)= \sum_{i=0:\infty}  g^{(i)}(s) ~ (t-s)^i/ i!  \simeq  g(s) +  g'(s) ~ (t-s) + {g''(s)}~ (t-s)^2/{2} ,
\ees
where by inserting the expansion in Equation (\ref{Er}) we have
\begin{align*}
Err(s)\leq \sum_{i=0:\infty} |g^{(i)}(s)| ~ Err^i(s)/ i!  \simeq  |g(s)| ~ Err^0(s)+  |g'(s)| ~ Err^1(s)+ \left| {g''(s)} \right| ~ Err^2(s)/{2}   ,
\end{align*}
when
\begin{align} \label{Ei}
Err^i(s)=\left| \int_{-1}^1     \ln(\|  t-s  \|) \, \| t - s \|^i \, dt - \sum_{k=1:n} \ln(\|  t_k-s  \|) \, \| t_k - s \|^i w_k \right| ~.
\end{align}
From Equation (\ref{Ei}), one can see $Err^i$ is the error of quadrature rule for function $f^i(r)=\ln(r) r^i$  from $-1$ to $1$ where $r=\| t-s \|$. We know $f^0(r)$ is unbounded at $r=0$ while $f^1(r)$ is continuous and $f^i(r)$ is differentiable at this point for $i \geq 2$. Since the Gaussian quadrature rule is more accurate for smoother functions, it seems the error of the rule decreases when $i$ increases. Graph of $Err^0, Err^1$ and $Err^2$ are shown in Figure \ref{fig2} for $n=16$ and $s\in[0,1]$. Since $Err^i$ is symmetric with respect to $s$, then $Err^i(s)=Err^i(-s)$ for $s<0$. The figure states $Err^0 \leq  Err^1 \leq Err^2$ almost everywhere.
Therefore, $Err^0$ is the main part of the error, $Err$, and the following semi-equality is valid
\be\label{bou}
Err(s) \simeq  |g(s)| ~ Err^0(s) .
\ee    

%\begin{figure}%[h]
%\begin{center}
%\includegraphics[width=8cm]{on_element.png}
%\caption{
%Computational domain $\Omega$ and its boundary, $\Gamma$. Source point, $\mbf{p}$, field point ,$\mbf{q}$, and distance vector, $\mbf{r}$ are shown in this picture.  Field points $\mbf{q}_{j+1}$ and $\mbf{q}_{j}$ are connected to each other by curved boundary element $\Gamma_j$.  
%}\label{fig0}
%\end{center}
%\end{figure}
%%%%%%%%%%%%%%%%%%%%%%%%%%%%%%%%%%%%%%%%%%%%%%%%%%%%
\subsection{Optimal value of $s$}\label{sec.op}
We are going to find the optimal value of $s$ minimizes the error function in this section, i.e.
\be\label{t0s}
s_{opt}=\arg \min_{s\in[-1,1]} \, Err(s) ~.
\ee
Thanks to semi-equality (\ref{bou}), the main source of the error is $Err^0$, and the error takes its minimum around zeros of $Err^0$. Consequently, the seeking optimal value also satisfies
\be\label{t0s2}
s_{opt} \simeq \arg \min_{s\in[-1,1]} \, Err^0(s)~.
\ee
Graph of $Err^0$ is depicted in Figure \ref{fig2} for  $s \in[0,1]$ and $n=16$. From Figure \ref{fig2} and this fact that $Err^0(s)$ is symmetric with respect to $s$, the zeros of $Err^0$ can be reported as
\begin{align}\label{con}
s = 
\left\{ \begin{array}{ll}
		%\left\{ \begin{array}{ll}
        %\pm0.123, \pm0.245, \pm0.474, \pm0.579, \pm0.761, & \\
         %\pm0.835, \pm0.942, \pm0.98, ~~~
         %\end{array} \right. 
        %& \mbox{if $n=8$} ,\\
        %& \\
        %\left\{ \begin{array}{ll} 
        \pm0.063, \pm0.128, \pm0.251, \pm0.312, \pm0.430, \pm0.486, & \\
        \pm0.593, \pm0.643, \pm0.735, \pm0.777, \pm0.850, \pm0.882, &  \\
         \pm0.935, \pm0.956, \pm0.985, \pm0.995 .    
        %\end{array} \right. 
        % & \mbox{if $n=16$} .
        \end{array} \right. 
\end{align}
These points satisfy Equation (\ref{t0s2}) and they are appropriate candidates for $s_{opt}$ when $n=16$. The efficiency of these points is checked numerically by the forthcoming example.% and one of them will be selected as the best one satisfying Equation (\ref{t0s}). Then a numerical example is presented in continue to find it. 

\begin{figure}%[h]
\begin{center}
\includegraphics[width=7.5cm]{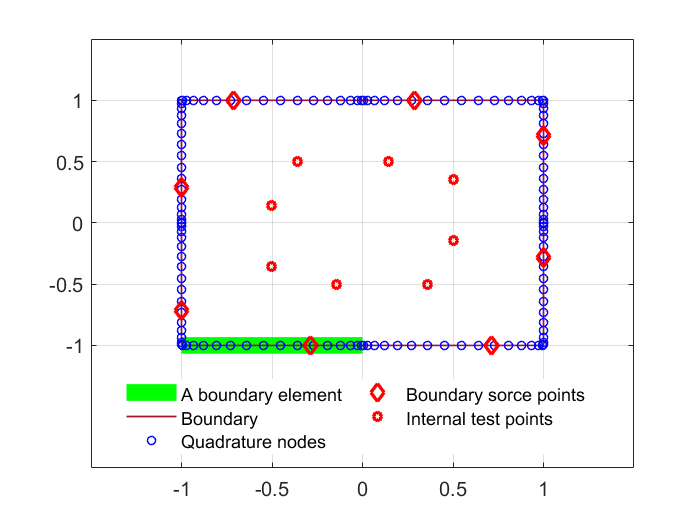} 
\includegraphics[width=7.5cm]{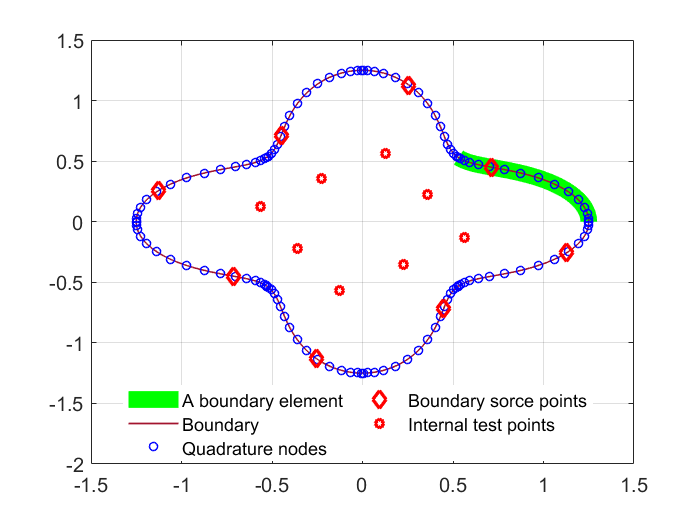} 
\caption{
Square (left) and flower-like (right) domains. The boundary elements are straight and curved for the domains, respectively.
}\label{fig_squ}
\end{center}
\end{figure}

\begin{figure}%[h]
\begin{center}
\includegraphics[width=13cm]{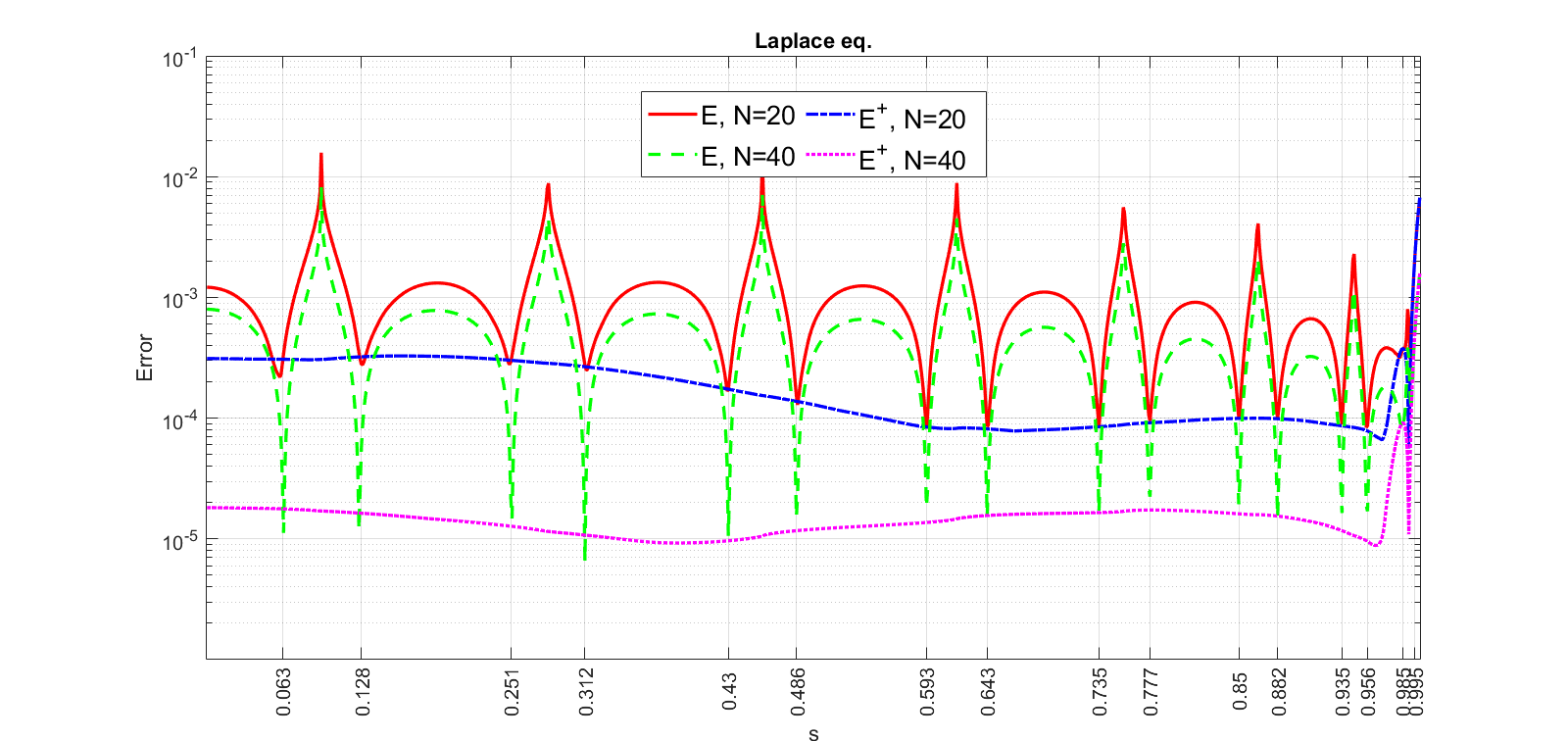}
\includegraphics[width=13cm]{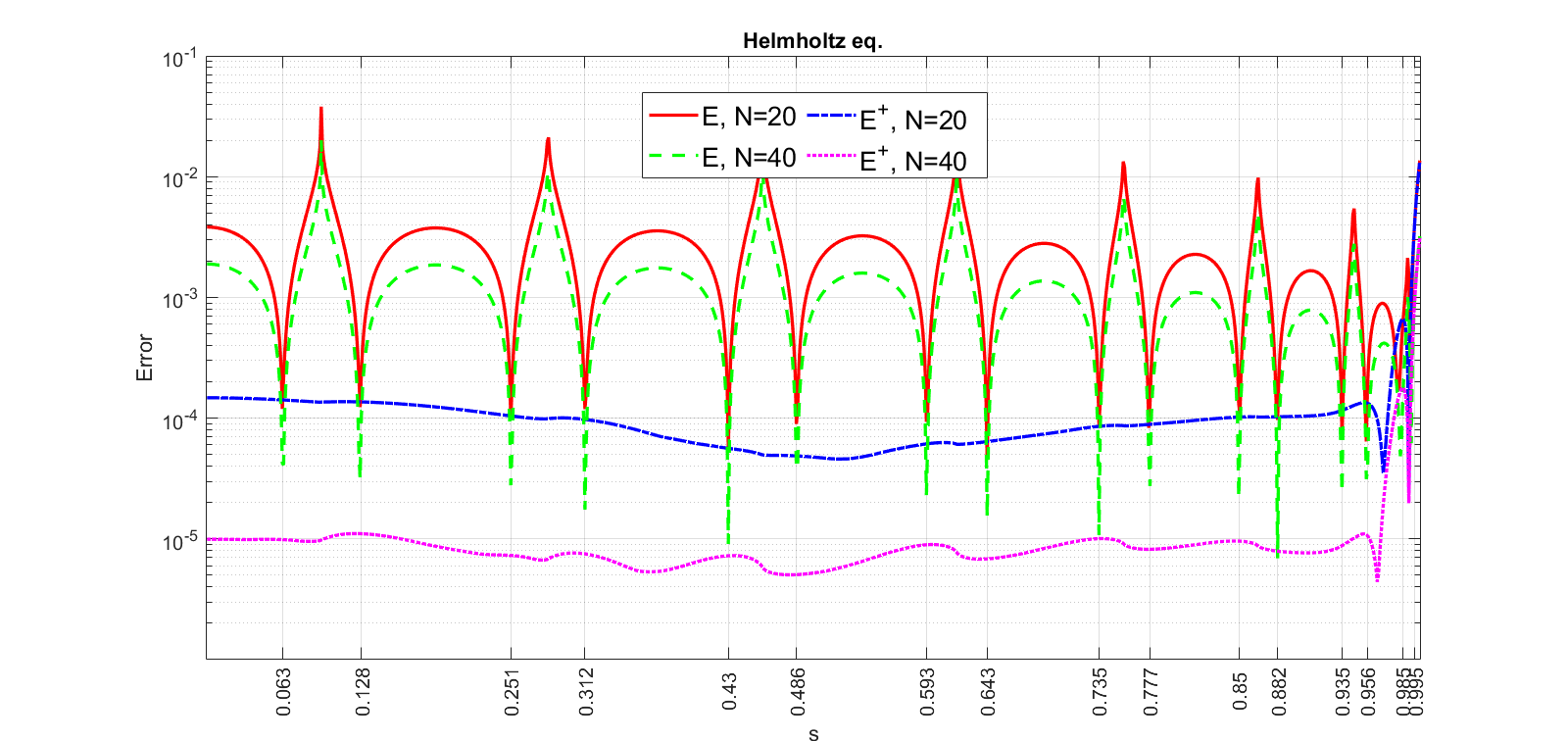}
\includegraphics[width=13cm]{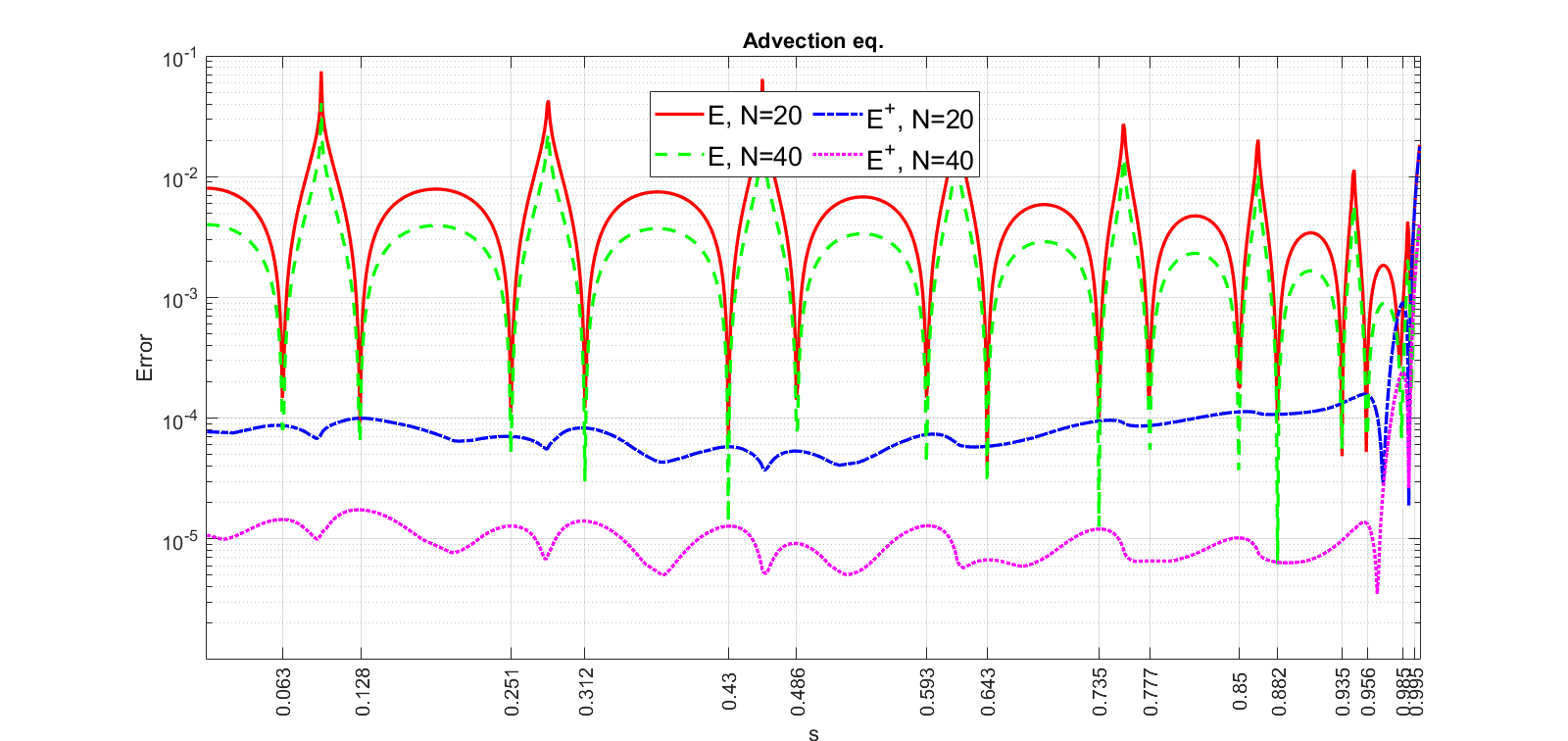}
\caption{
Graph of error of the radial BEM to solve Laplace (upper), Helmholtz (middle) and Advection (lower) equations in field of $s$ for the square domain. The error is reported when singular integrals are avoided ($E$) and they are calculated analytically ($E^+$). It can be seen the error is reduced significantly at roots of $Err^0$ when the singular integrals are neglected. However, the error doesn't change significantly for $s\in[0,0.9]$ when the singular integrals are calculated analytically. 
}\label{fig333}
\end{center}
\end{figure}

\begin{figure}%[h]
\begin{center}
\includegraphics[width=13cm]{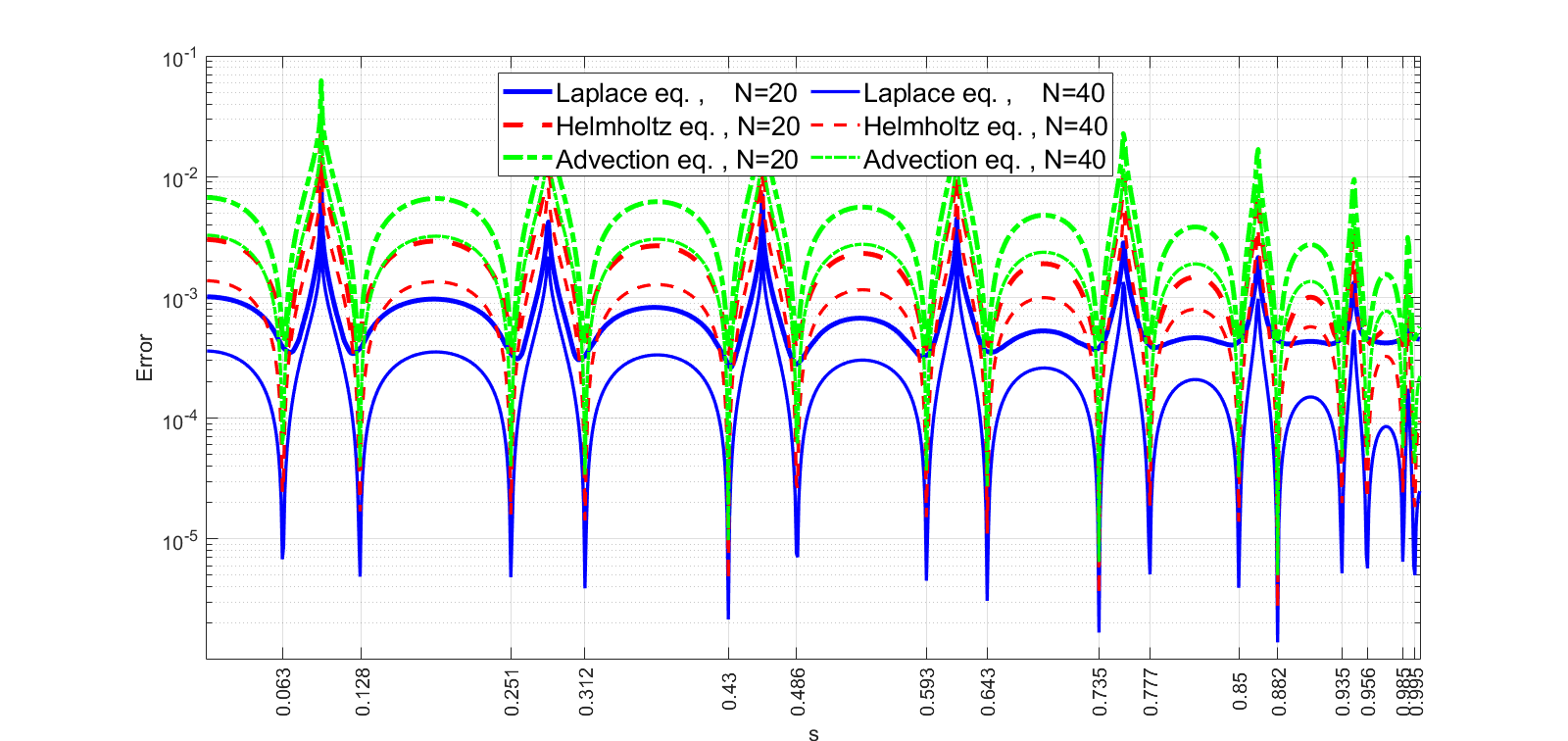}
\caption{
Graph of error of the radial BEM to solve Laplace, Helmholtz and Advection equations in field of $s$ for the flower-like domain. The error is reported when singular integrals are avoided. It can be seen the error is reduced significantly at roots of $Err^0$. 
}
\label{fig444}
\end{center}
\end{figure}

%%%%%%%%%%%%%%%%%%%%%%%%%%%%%%%%
%\begin{Example}\label{ex1}
Let $\Omega=[-1,1]^2$ and its boundary is discretized to $N$ equal boundary elements. Then, $N$ source points 
\be\label{col}
\mbf{p}_1, \mbf{p}_2 , \mbf{p}_3,  ... , \mbf{p}_{N}  ,
\ee     
are selected on the boundary elements where $ \mbf{p}_{j} \in \Gamma_j$ and $\mbf{p}_{j}=q_j(s)$ for a $s\in[-1,1]$. PDEs
\begin{align}
\hbox{Laplace eq. :} & ~~~ \nabla^2 u  = 0 , \label{pde1}\\
\hbox{Helmholtz eq. :} & ~~~ \nabla^2 u  - 2 u  =0 , \label{pde2}\\
\hbox{Advection eq. :} & ~~~ \nabla^2 u  +2 u_x - u_y - 5 u =0 , \label{pde3} 
\end{align}
are considered with Dirichlet boundary condition when the exact solution of the PDEs are
\begin{align}
u(x,y) &=\exp(x)\cos(y)   , \label{anal1}\\
u(x,y) &=\exp(x+y)   , \label{anal2}\\
u(x,y) &=\exp(x -y)   , \label{anal3}
\end{align}
respectively. Parameters $x$ and $y$ refer to the components of $\mathbb{R}^2$.
So the radial BEM is applied and the flux function, $v$, is obtained at the boundary source points. After that, the potential function is obtained numerically at internal test points 
\be\label{intest}
\hat{\mbf{p}}_{j}=\frac{1}{2} \mbf{p}_j ,
\ee
for $j=1, 2, ..., N$ by Equation (\ref{eq2.7}). The boundary source points and the internal test points are shown in Figure  \ref{fig_squ} (left) for $N=8$ when $s=0.430$. Boundary integrals of the method are calculated by GQR with $n=16$ quadrature nodes, and error of the method is calculated by error function
\be\label{err}
E(s)=\frac{1}{N} \sum_{j=1:N} |\hat{u}[j]-u(\hat{\mbf{p}}_j)| ,
\ee
where $u$ is the exact solution, and $\hat{u}[j]$ is the numerical solution at $\hat{\mbf{p}_j}$. 
Graph of the error function, $E$, is shown in Figure \ref{fig333} in field of $s\in[0,1]$ when $N=20$ and $40$. Note that, the error function is symmetric with respect to $s$ and consequently $E(s)=E(-s)$ for $s\in[-1,0]$. 
Gaussian RBFs 
\begin{align}\label{rbf}
\begin{matrix}
\phi_{j}(\mbf{q}) &= \exp(-\epsilon^2 \| \mbf{q} - \mbf{p}_{j} \|^2) ,
\end{matrix}
\end{align}
are applied in this example to approximate boundary values for $j=1, 2, ...,N$ where the shape parameter of the RBFs is evaluated as
\be\label{eps}
\epsilon =N/\sqrt{1000} .
\ee
Numerical experiments, presented in Figure \ref{fig333}, show error function $E$ is minimized locally at the roots of $Err^0$ when $s$ is not close to $1$. The error function also is calculated when singular integrals, including $\ln(r)$, are calculated analytically by the idea presented in \cite{hosseinzadeh2014simple}. By obtaining the singular integrals analytically, $Err^0$ is vanished from error formula (\ref{Er}) and $Err^1$ will be the main part of the error of the quadrature rule, $Err$. The modified error function is shown in Figure \ref{fig333} by notation $E^+$. From Figure \ref{fig333}, $E^+(s)$ does not vary significantly for $s\in[0,0.9]$. Therefore, those roots of $Err^0$ absolutely smaller than $0.9$ can be supposed as optimal values of $s$ for the square domain. The experiment also is done for flower-like region
\begin{align} \label{flower}
\Omega=\{ (r \cos(\theta), r\sin(\theta))\in\mathbb{R}^2 ~| ~ 0\leq \theta \leq 2\pi , ~0 \leq r \leq (1+0.25\cos(4 \theta)) \},
\end{align}
when its boundary is discretized to $N$ semi-equal curved boundary elements. The region is shown in Figure \ref{fig_squ} (right) for $N=8, n=16$ and $s=0.430$. Graph of error function $E$ is shown in Figure \ref{fig444} in field of $s\in[0,1]$ for PDEs (\ref{pde1})-(\ref{pde3}) with Dirichlet boundary condition imposed on the boundary source points by exact solutions (\ref{anal1})-(\ref{anal3}). Figure \ref{fig444} verifies that $E$ is minimized at zeros of $Err^0$ also for domains with curved boundaries.
Then, from Equation (\ref{con}) and the above statements the following corollary is valid.

\begin{cor}\label{cor1}
The error of radial BEM is minimized at $s= 0.430$ when the Gaussian quadrature rule with $n=16$ quadrature nodes is applied for numerical integration. In this situation, we have $s_{opt}= 0.430$.
\end{cor}

Corollary \ref{cor1} proposes an optimal value for $s$ when $n=16$. But, the idea can be extended for the other values of $n$, and find
\begin{align}\label{t0}
s_{opt}= 
\left\{ \begin{array}{ll}
        \pm 0.453 ,~~~& \mbox{if $n=4$} ,\\
        \pm 0.474 ,~~~& \mbox{if $n=8$} , \\
        \pm 0.430 ,~~~& \mbox{if $n=16$} ,\\
        \pm 0.521 ,~~~& \mbox{if $n=32$} ,\\
        \pm 0.410 ,~~~& \mbox{if $n=64$} .
        \end{array} \right. 
\end{align}
It seems larger values of $n$ lead to more accurate results, but it is not valid, generally. The next subsection is devoted to this subject.

%%%%%%%%%%%%%%%%%%%%%%%%%%%%%%%%%%%%%%%%%%%%%%%%
\subsection{Optimal value of $n$}\label{sec.op2}
\begin{figure}%[h]
\begin{center}
\includegraphics[width=7.5cm]{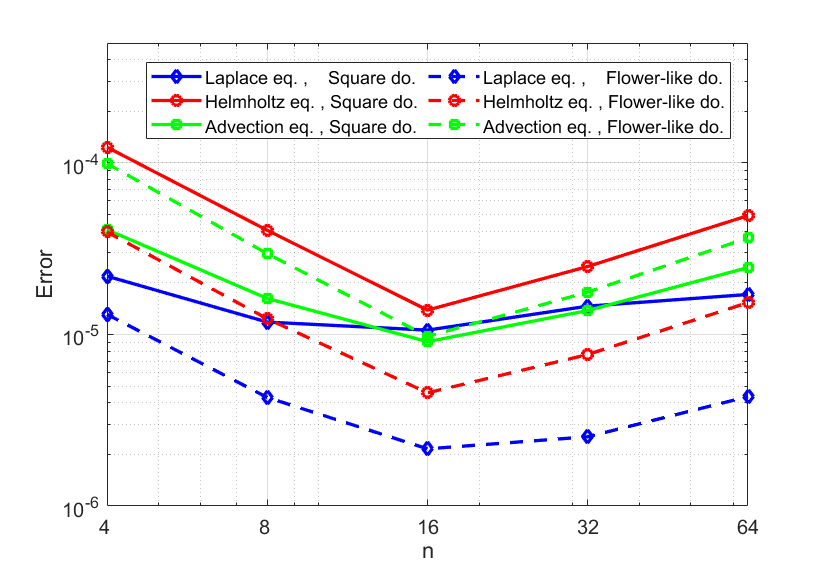}
\includegraphics[width=7.5cm]{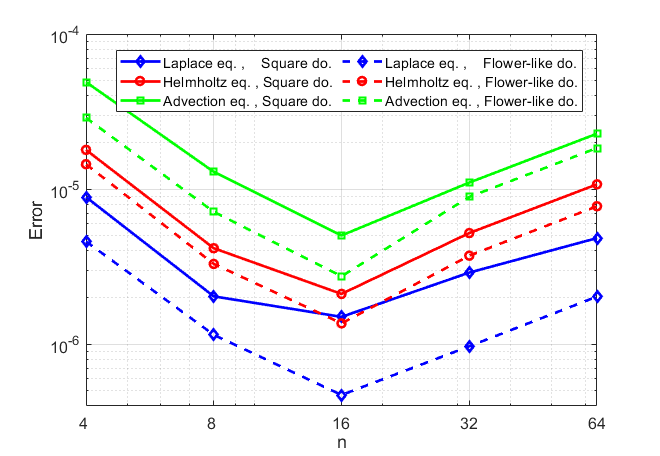}
\caption{
The error of the radial BEM in field of $n$ for $N=40$ (left) and $N=80$ (right). From this figure, one can see the error is minimized at $n=16$.
}\label{fig555}
\end{center}
\end{figure}

Numerical experiments show the quadrature rule is not sufficiently accurate for $n \leq 8$. It is happened because semi-equality (\ref{bou}) is valid only for sufficiently large values of $n$. At the other hand, $s_{opt}$ is very sensitive to the number of the digits for $n\geq 32$, and the proposed optimal values, presented in Equation (\ref{t0}), do not enhance the accuracy, numerically. The next example shows the best value of $n$ is nearly $16$ when $N$ is sufficiently large. Let PDEs (\ref{pde1})-(\ref{pde3}) are valid in the square and flower-like domains when Dirichlet boundary condition is imposed on their boundary by exact solutions (\ref{anal1})-(\ref{anal3}).
The boundary of computational domains is discretized to $N$ semi-equal boundary elements and Gaussian RBFs (\ref{rbf}) with shape parameter (\ref{eps}) are applied to approximate boundary values. The PDEs are solved by the radial BEM when the optimized boundary source points $\mbf{p}_j=q_j(s_{opt})$ are applied for $j=1, 2, ..., N$. The error of the radial BEM at internal test points (\ref{intest}) is calculated by formulation (\ref{err}) and is depicted in Figure \ref{fig555} for $n=4, 8, 16, 32$ and $64$ when $N=40$ and $80$. From this figure, the error is minimized at $n=16$ for the case studies. Therefore, larger values of $n$ do not enhance the accuracy of BEM and one can set $n=16$ and increase $N$ to get more accurate results. %This fact also is valid for larger values $N$ which shows there is no need to $n>16$ when $N$ is sufficiently large.

\section{Numerical experiments}\label{Sec4}
From Section \ref{Sec3}, the radial BEM is accurate if the boundary integrals are calculated by the GQR with $n=16$ quadrature nodes. In this section, the radial BEM is investigated, by several numerical examples. At first, the precision of radial BEM is studied in Subsection \ref{rbfs-rbem} for some well-known RBFs. Then, the radial BEM is compared with the standard BEM in terms of accuracy in Subsection \ref{r-l-bem}. And it is compared with RBF collocation method  in terms of stability in Subsection \ref{r-r-bem}. Numerical results are also analysed mathematically to highlight the advantages of the new method. %The shape parameter of the RBFs is considered as (\ref{eps}) in this section. %In this example Gaussian RBFs are applied for the approximation.
%In the radial BEM, boundaries of computational domains are discretized to $N$ semi-uniform boundary elements and some boundary source points, respect to Corollary (\ref{cor1}), are selected on them. These source points minimize error of numerical integration and enrich BEM very much. %More details are also same as what are described in the previous section. 
%After that, numerical solution of BEM for advection-diffusion equation is studied in subsection \ref{Ex3}. Accuracy of BEM is clarified there when RBFs and Gaussian quadrature rule are applied to approximate boundary values and boundary integrals, respectively.

%%%%%%%%%%%%%%%%%%%%%%%%%%%
\subsection{Radial BEM for several RBFs}\label{rbfs-rbem}
Let square domain $\Omega=[-1,1]^2$ be considered as computational domain and its boundary is discretized to $N$ straight boundary elements. Figure \ref{fig_squ} (left) shows the boundary when it is discretized to $N=8$ boundary elements and $n=16$ quadrature nodes are selected on each boundary element to calculate the boundary integrals. Consequently, Laplace equation (\ref{pde1})
is solved by radial BEM and error of the method is calculated by formulation (\ref{err}) for internal test points (\ref{intest}). The error is reported in Table \ref{tab2} for Gaussian, inverse multiquadric (IMQ), thin plate spline (TPS), poly-harmonic spline (PHS) and a local radial basis functions \cite{wendland2004scattered, fasshauer2007meshfree}.
The first two RBFs are infinitely smooth functions containing shape parameter $\epsilon$ evaluated as (\ref{eps}). To study the accuracy of radial BEM, Dirichlet and mixed boundary conditions are imposed at the boundary source points by exact solution (\ref{anal1}).
In the Dirichlet boundary condition, the potential value is known at the boundary source points. But, in the mixed boundary condition, the potential value is available only at the upper and lower sides of the square domain while flux is known at the left and right sides. From Table \ref{tab2}, the error of radial BEM is smaller than $10^{-5}$ for Gaussian, IMQ and PHS RBFs when $N \geq 64$ and Dirichlet boundary condition is imposed. Also, the error is smaller than $10^{-3}$ for the RBFs when $N \geq 64$ and the mixed boundary condition is applied. Therefore, radial BEM is accurate for the square domain.

To check the accuracy of the radial BEM for domains with curved boundaries, flower-like region (\ref{flower}) is considered as a computational domain. The boundary of the domain is discretized to $N$ curved semi-equal boundary elements as is shown in Figure \ref{fig_squ} (right) for $N=8$. Similar to the square domain, $N$ optimal boundary source points are selected on the boundary elements and boundary integrals are calculated by the GQR with $n=16$ nodes. The error of the radial BEM is calculated by equation (\ref{err}) for PDE (\ref{pde1}) and it is reported in Table \ref{tab3} for Dirichlet and mixed boundary conditions. In the mixed boundary condition potential and flux functions are known for the upper and lower half boundary points, respectively. From this table, Gaussian, IMQ and PHS RBFs yield accurate results when $N\geq 64$. 

Regarding tables \ref{tab2} and \ref{tab3}, radial BEM is able to solve the PDE accurately when infinitely smooth RBFs are applied and the boundary integrals are calculated by the GQR with $n=16$ nodes. Note that, radial BEM is more accurate for flower-like domain versus square domain because it does not have any corners and BEM works better for regions with smooth boundaries \cite{sedaghatjoo2013use}. % while the optimal boundary points are applied as source points.
It must be mentioned that since the approximation power of smooth RBFs is significantly higher than piecewise polynomials, radial BEM is significantly more accurate than standard BEM. As a special case, the linear BEM can be supposed as a radial BEM when local RBF $\phi(r)=(1-r)_+$ is applied. In fact, the local RBF works similar to a linear function on the one-dimensional boundary elements. Therefore, from tables \ref{tab2} and \ref{tab3} radial BEM can be significantly more accurate than linear BEM when RBFs with high approximation power are utilized. The next subsection studies this subject with more details.

 %  are able to enhance the accuracy of BEM significantly, specially when the boundary integrals are calculated accurately.

\begin{table}[ht]
  \begin{center}
    \caption{Numerical results of radial BEM for the square domain. }\label{tab2}
    \begin{tabular}{cccccccccc} % <-- Changed to S here.
      \hline
      \hline
           &     &     &      & $N$  &      &   \\
      \cline{3-7}
      B.C.  & RBF  &  $16$ & $32$ & $64$ & $128$ & $256$  \\
      \hline
            & Gaussian  & $3.73e{-4}$ & $2.08e{-5}$ & $2.73e{-6}$ & $7.29e{-7}$ &$5.25e{-7}$ \\
            & IMQ         & $4.92e{-4}$ & $3.06e{-5}$ & $3.98e{-6}$ & $7.42e{-7}$ &$2.54e{-7}$ \\
      Dirichlet & TPS  & $7.80e{-3}$ & $4.54e{-4}$ & $4.15e{-5}$ & $4.40e{-6}$ &$5.09e{-7}$ \\
             & PHS  & $1.20e{-3}$ & $6.89e{-5}$ & $6.63e{-6}$ & $8.82e{-7}$ &$2.53e{-7}$ \\
%              & Local $C^0$  &  $3.00e{-2}$ & $9.40e{-3}$ & $2.70e{-3}$ & $7.04e{-4}$ &$1.81e{-4}$ \\
              & Local $C^0$  &  $1.29e{-2}$ & $3.40e{-3}$ & $7.45e{-4}$ & $1.86e{-4}$ &$4.41e{-5}$ \\
                    \hline
            & Gaussian  &  $1.10e{-2}$ & $2.70e{-3}$ & $7.16e{-4}$ & $2.02e{-4}$ &$5.00e{-5}$ \\
            & IMQ  &  $1.17e{-2}$ & $3.10e{-3}$ & $8.73e{-4}$ & $2.48e{-4}$ &$7.22e{-5}$ \\
      Mixed  & TPS  &  $1.34e{-2}$ & $2.80e{-3}$ & $8.30e{-4}$ & $2.45e{-4}$ &$7.27e{-5}$ \\
            & PHS  &  $1.15e{-2}$ & $3.30e{-3}$ & $9.31e{-4}$ & $2.66e{-4}$ &$7.70e{-5}$ \\
%            & Local $C^0$  &  $4.32e{-2}$ & $1.72e{-2}$ & $6.00e{-3}$ & $2.00e{-3}$ &$6.19e{-4}$ \\
            & Local $C^0$  &  $4.67e{-2}$ & $2.09e{-2}$ & $9.90e{-3}$ & $4.90e{-3}$ &$2.50e{-3}$ \\
                    \hline
                    \hline
   \end{tabular}
  \end{center}
\end{table}

\begin{table}[ht]
  \begin{center}
    \caption{Numerical results of radial BEM for the flower-like domain. }\label{tab3}
    \begin{tabular}{cccccccccc} % <-- Changed to S here.
      \hline
      \hline
           &     &     &      & $N$  &      &   \\
      \cline{3-7}
      B.C.  & RBF  & $16$ & $32$ & $64$ & $128$ & $256$ \\
      \hline
            & Gaussian   & $1.30e{-3}$ & $1.48e{-5}$ & $5.89e{-7}$ & $2.40e{-7}$ &$1.07 e{-7}$ \\
            & IMQ          & $1.20e{-3}$ & $1.52e{-5}$ & $5.88e{-7}$ & $2.34e{-7}$ & $1.06e{-7}$  \\
      Dirichlet & TPS  & $2.70e{-3}$ & $1.89e{-4}$ & $2.21e{-5}$ & $2.67e{-6}$ & $3.07e{-7}$  \\
             & PHS  & $1.30e{-3}$ & $1.96e{-5}$ & $7.28e{-7}$ & $2.32e{-7}$ & $1.07e{-7}$  \\
%              & Local $C^0$  &  $2.63e{-2}$ & $7.30e{-3}$ & $1.90e{-3}$ & $4.76e{-4}$ & $1.19e{-4}$  \\
              & Local $C^0$  &  $1.80e{-2}$ & $3.30e{-3}$ & $1.10e{-3}$ & $1.35e{-4}$ & $4.90e{-5}$  \\
                    \hline
            & Gaussian  &  $8.30e{-3}$ & $1.42e{-4}$ & $1.26e{-6}$ & $2.49e{-7}$ & $2.31e{-7}$  \\
            & IMQ  &         $3.50e{-3}$ & $1.68e{-4}$ & $9.15e{-7}$ & $2.51e{-7}$ & $1.49e{-7}$  \\
      Mixed  & TPS  &    $1.01e{-2}$ & $6.84e{-4}$ & $9.86e{-5}$ & $1.19e{-5}$ & $1.61e{-6}$  \\
            & PHS  &        $1.07e{-2}$ & $1.92e{-4}$ & $7.90e{-6}$ & $5.12e{-7}$ & $1.74e{-7}$  \\
%            & Local $C^0$  &  $2.66e{-2}$ & $4.80e{-3}$ & $1.20e{-3}$ & $3.09e{-4}$ & $7.72e{-5}$  \\
            & Local $C^0$  &  $1.21e{-2}$ & $4.00e{-3}$ & $1.40e{-3}$ & $7.30e{-4}$ & $4.37e{-5}$  \\
                    \hline
                    \hline
   \end{tabular}
  \end{center}
\end{table}

%%%%%%%%%%%%%%%%%%%%%%%%%%%%%%%
\subsection{Radial BEM versus linear BEM}\label{r-l-bem}
\begin{figure}%[h]
\begin{center}
\includegraphics[width=13cm]{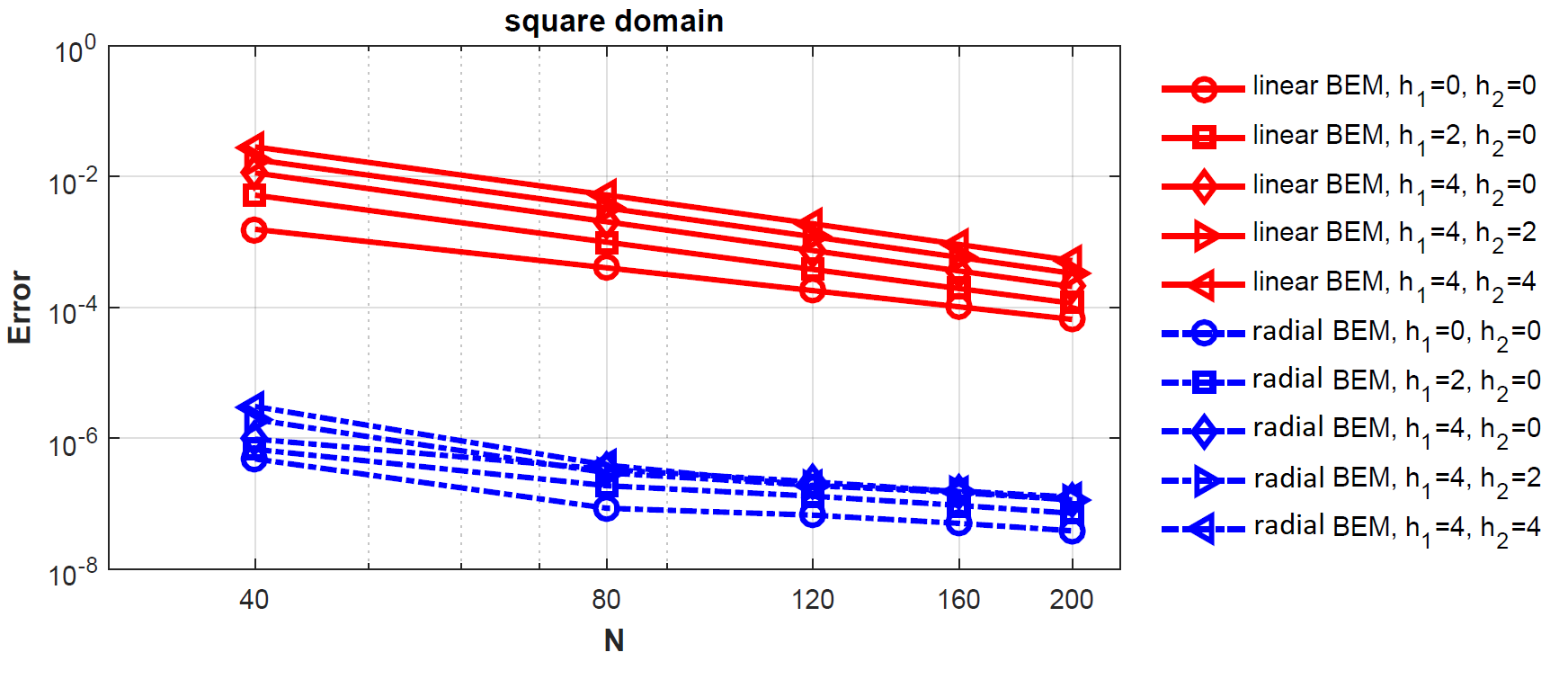} 
\includegraphics[width=13cm]{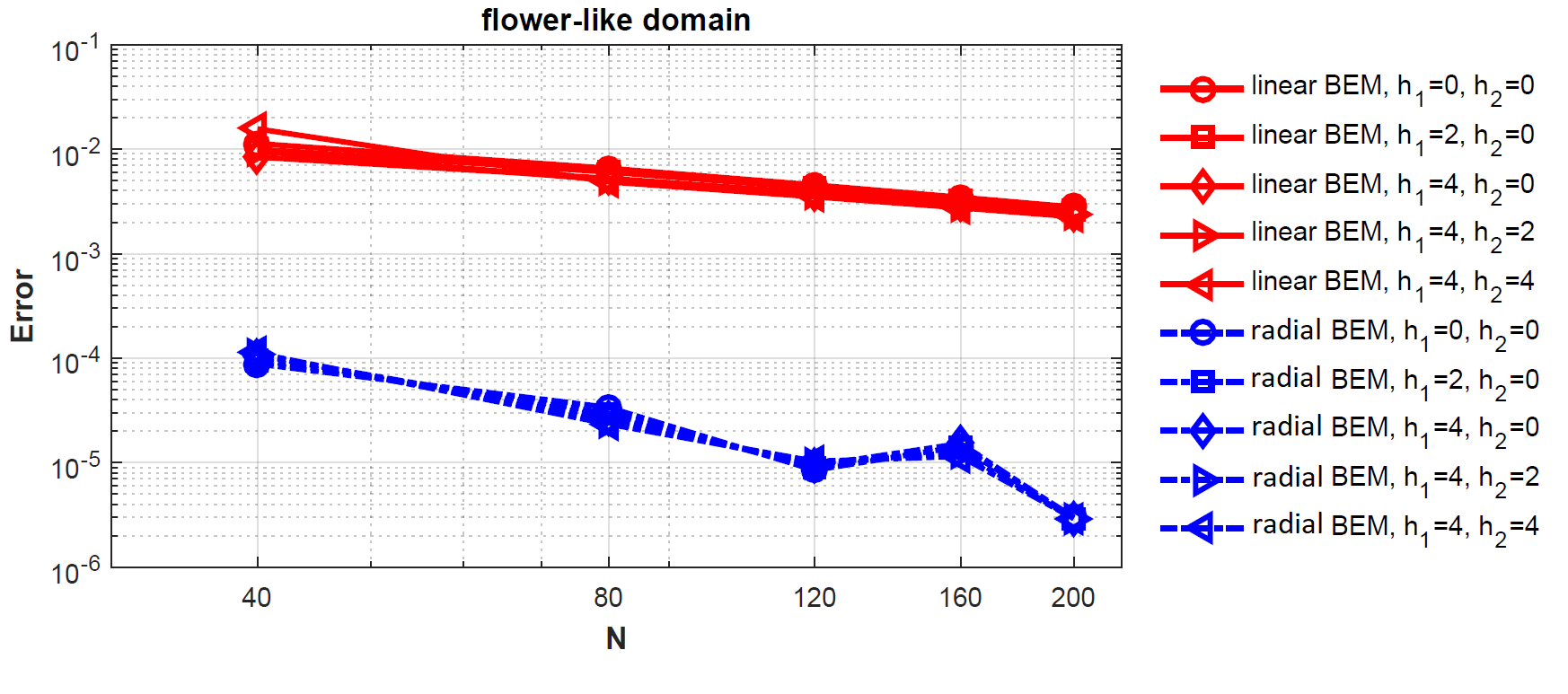} 
\caption{
Graph of error of the linear and the radial BEMs to solve PDE (\ref{PDE}) for several values of $h_1$ and $h_2$. Computational domain is considered as square (up) and flower-like (dawn) regions, and boundary integrals of the methods are calculated by Gaussian quadrature rule with $n=16$ quadrature nodes. It is clear that the error of the radial BEM is significantly smaller than the error of the linear BEM in the case studies.
}\label{fig6}
\end{center}
\end{figure}

In this subsection, the accuracy of the radial BEM is compared with that of the linear BEM for solving PDE (\ref{PDE}).
The square and flower-like regions considered as computational domains are shown in Figure \ref{fig_squ}. The boundary of the domains is discretized to $N$ straight boundary elements in the methods to compare the approximation power of smooth RBFs with that of linear polynomials when they are applied in BEM formulation.
The applied RBFs implemented in the radial BEM are
\begin{align}\label{rbf2}
\begin{matrix}
\phi_{2j-1}(\mbf{q}) &= \exp(-\epsilon^2 \| \mbf{q} - q_j(-s_{opt}) \|^2) , \\
\phi_{2j}(\mbf{q}) &= \exp(-\epsilon^2\| \mbf{q} - q_j(+s_{opt}) \|^2) , 
\end{matrix}
\end{align}
for $j=1, 2, ..., N$, while linear functions
\begin{align}\label{linear}
\begin{matrix}
\phi_{2j-1}(\mbf{q})  &= -\frac{1}{2 \, s_{opt}}(t-s_{opt}) , & \\
\phi_{2j} (\mbf{q})  &= + \frac{1}{2 \, s_{opt}}( t+s_{opt}) . &
\end{matrix}
\end{align}
are implemented in the linear BEM. Note that, two boundary source points are considered on each boundary element as $\mbf{p}_{2j-1}= q_j(-s_{opt})$ and $\mbf{p}_{2j}= q_j(s_{opt})$. Thanks to Corollary \ref{cor1}, we set $s_{opt}=0.430$ to calculate the boundary integrals accurately by the GQR with $n=16$ quadrature nodes. The shape parameter of the RBFs is evaluated as (\ref{eps}).
Dirichlet boundary condition is imposed on the boundary points by exact solution (\ref{anal2}),
and the mean value of error, Equation (\ref{err}), is calculated for the methods at internal points (\ref{intest}). The error is shown in Figure \ref{fig6} when the number of boundary elements varies from $40$ to $200$.
It can be found from the figure that the error of the radial BEM is significantly less than the error of the linear BEM which shows higher approximation power of the RBFs versus the linear polynomials.

%%%%%%%%%%%%%%%%%%%%%%%%%%%%%%%%%%%
\subsection{Radial BEM versus RBF collocation method}\label{r-r-bem}
\begin{figure}%[h]
\begin{center}
\includegraphics[width=7.5cm]{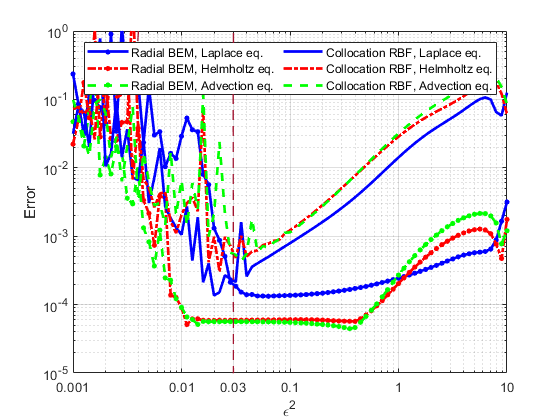} 
\includegraphics[width=7.5cm]{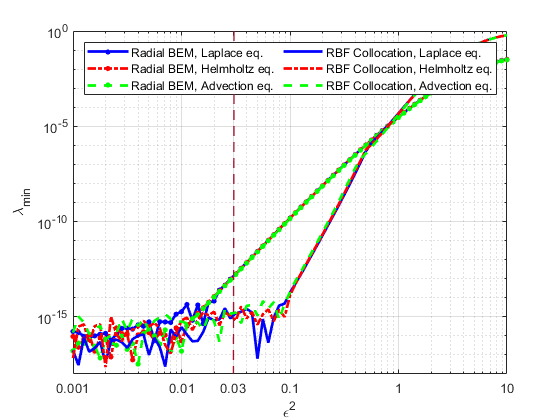} \\
\includegraphics[width=7.5cm]{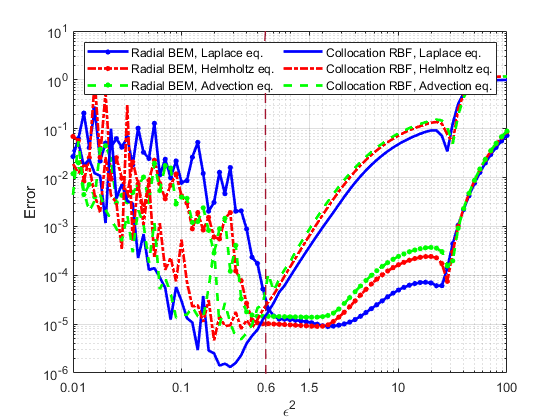} 
\includegraphics[width=7.5cm]{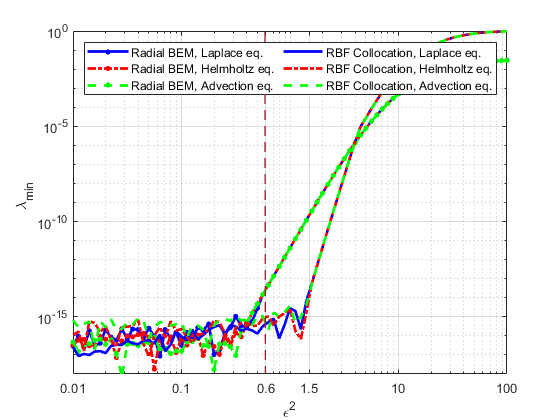} \\
\includegraphics[width=7.5cm]{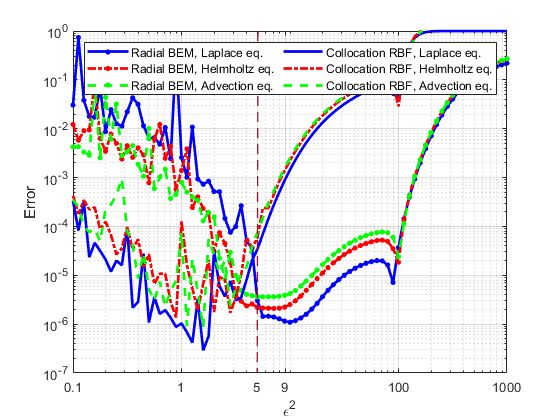} 
\includegraphics[width=7.5cm]{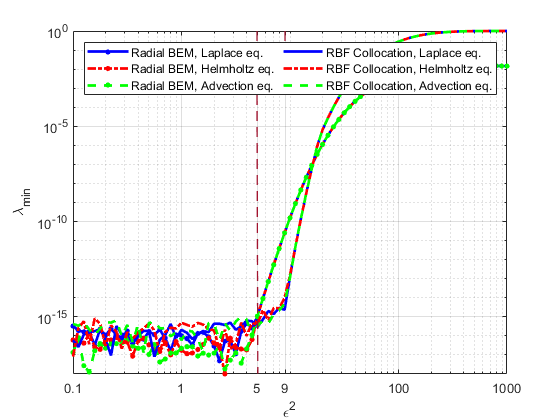} 
\caption{
Graph of the error of the radial BEM and the RBF collocation method to solve PDEs  (\ref{pde1})-(\ref{pde2}) for $h_{min}=0.2, 0.1$ and $0.05$ from up to down, respectively. It is clear that the results of the radial BEM is significantly more stable than the results of the RBF collocation method. 
}\label{fig666}
\end{center}
\end{figure}

%%%%%%%%%%%%%%%%
%\begin{Example}\label{ex4}
In this subsection accuracy of radial BEM is compared with that of RBF collocation method to solve PDEs (\ref{pde1})-(\ref{pde2}).
The square domain, $\Omega=[-1,1]^2$, is considered as the computational domain, and the Dirichlet boundary condition is imposed on its boundary by the exact solutions (\ref{anal1})-(\ref{anal3}).
The boundary of the domain is discretized to $N$ straight boundary elements in the radial BEM and $N$ boundary source points are selected on the boundary elements as $\mbf{p}_j=q_j(s)$ for $j=1, 2, ..., N$ where $s=0.430$ as is mentioned in Corollary \ref{cor1}. Then, boundary values are approximated by Gaussian RBFs (\ref{rbf}), and boundary integrals are calculated by the GQR with $n=16$ quadrature nodes.
In RBF collocation method, $N^2/16$ center points are selected in $\Omega$ and on its boundary, uniformly, and numerical solution
\[
\bar{u}(\mbf{x})=\sum_{j=1:N^2/16} \gamma_j \phi(r_j(\mbf{x})) ,
\]
is found such that it satisfies the PDE and the Dirichlet boundary condition at internal and boundary center points, respectively \cite{vsarler2007global, vsarler2012radial}. This method also is called direct collocation and Kansa's methods in the literature \cite{koupaei2018finding, rosenfeld2019mesh}. Gaussian RBFs are applied here for approximation. It is well-known that the RBF collocation method is very sensitive to the shape parameter of the RBFs and its best results are obtained at that shape parameters spoil its stability. Errors of the radial BEM and the RBF collocation method at internal test points (\ref{intest}) are calculated by Equation (\ref{err}), and they are depicted in Figure \ref{fig666} in field of $\epsilon^2$. The results are presented for $N=20, 40$ and $80$ where $h_{min}=0.2, 0.1$ and $0.05$ for them, respectively. From this figure, numerical results of the radial BEM are significantly more accurate than the RBF collocation method when $\epsilon^2\geq 0.03, 0.6$ and $5$ for $h_{min}=0.2, 0.1$ and $0.05$, respectively.
Note that, the reported shape parameters are lower bounds for ill-conditioning of system of linear equations (\ref{final}) and the radial BEM is not stable for shape parameters smaller than those. The lower bound for the absolute value of eigenvalues of system (\ref{final}), denoted by $\lambda_{min}$, is presented in Figure \ref{fig666} (right). From Figure \ref{fig666}, $\lambda_{min} \leq 10^{-13}$ for $\epsilon^2\leq 0.03, 0.6$ and $5$ when $h_{min}=0.2, 0.1$ and $0.05$, respectively. The value of $\lambda_{min}$ is also presented in Figure \ref{fig666} (right) for the final system of the RBF collocation method.
From this figure, the system of the RBF collocation method is spoiled faster than the system of radial BEM when $\epsilon^2$ tends to $0$. It has happened because in radial BEM only boundary values are approximated by the RBFs but in the RBF collocation method domain values are approximated by them. Therefore, the dimension of the problem is reduced by one for the radial BEM and thanks to Table \ref{tab1} it is more stable than the RBF collocation method. This is the main advantage of the radial BEM versus the RBF collocation method. 
%Note that, boundary of the computational domain is one dimension and RBFs in Radial BEM approximate boundary values. Then RBFs in Radial BEM can be supposed as one dimensional functions while they are two dimensional functions for collocation RBF method. As a result, From Table \ref{tab1}, Gaussian RBFs in Radial BEM are more stable than those in the collocation method specially when $h_{min}$ tends to zero. Figure \ref{fig666} verifies this fact, numerically.

%%%%%%%%%%%%%%%%%%%%%%%%%%%%
%%%%%%%%%%%%%%%%%%%%%%%%%%%%
\section{Conclusion}\label{Sec5}
Radial BEM was proposed in this paper wherein radial basis functions (RBFs) were applied to approximate boundary values.
To overcome the singularity problem of the proposed method, a new distribution of boundary source points was proposed. It allowed us to use the standard Gaussian quadrature rule (GQR) for boundary integrals of BEM, safely.
The radial BEM was studied numerically for several RBFs including Gaussian RBFs, and it was compared with standard BEM and RBF collocation method via several examples.
The results show that the radial BEM is much more accurate than the standard BEM because of the high approximation power of infinitely smooth RBFs.
It is also significantly more stable than the RBF collocation method because of the dimension reduction of BEM.
The radial BEM was applied here for two dimensional linear partial differential equations (PDEs), but it is applicable for solving more complicated or higher dimensional ones.
The new approach is also applicable to other basis functions, such as sigmoid and inverse hyperbolic functions, which have been applied in machine learning problems, extensively.
Moreover, the new technique applied removing the singularity of the boundary integrals is suitable to handle singularity of other boundary integral equations concluded singular kernels.
%%%%%%%%%%%%%%%%%%%%%%%%%%%%%%%%%%%%%%%%%%%%%%%%%%%%%
%%%%%%%%%%%%%%%%%%%%%%%%%%%%%%%%%%%%%%%%%%%%%%%%55
%\begin{small}
%\begin{thebibliography}{99}
%\section*{References}
\bibliographystyle{elsarticle-num}
\bibliography{mybibfile}

%\bibitem{Atk1997} K. E. Atkinson, The Numerical Solution of Integral Equations of the Second Kind, Cambridge Monographs on Applied and Computational Mathematics, Cambridge, New York, 1997.

%\end{thebibliography}
%\end{small}
%%%%%%%%%%%%%%%%%%%%%%

%%%%%%%%%%%%%%%%%%%%%%
\end{document}